\documentclass[leqno,11pt]{article}
\usepackage{amsmath,amssymb,amscd,amsthm,amsfonts,epsfig}

\swapnumbers

\newtheorem{thm}[equation]{Theorem}
\newtheorem{prop}[equation]{Proposition}
\newtheorem{lem}[equation]{Lemma}

\newtheorem{cor}[equation]{Corollary}
\newtheorem{rem}[equation]{Remark}
\newtheorem{df}[equation]{Definition}
\newtheorem{ex}[equation]{Example}

\renewcommand{\phi}{\varphi}
\renewcommand{\epsilon}{\varepsilon}

\newcommand{\BB}{\mathbb}
\newcommand{\g}{\mathfrak}
\newcommand{\pf}{\noindent {\it Proof. }}
\newcommand{\separate}{\vskip7pt}
\newcommand{\tr}{\operatorname{tr}}
\newcommand{\End}{\operatorname{End}}
\newcommand{\Hom}{\operatorname{Hom}}

\newcommand{\Ad}{\operatorname{Ad}}
\newcommand{\HC}{\operatorname{HC}}

\numberwithin{equation}{subsection}

\begin{document}

\title{\bf Geometric Methods in Representation Theory}
\author{Wilfried Schmid\thanks{Supported in part by NSF grant DMS-0070714}  \\ Lecture Notes Taken by Matvei Libine}
\maketitle

\begin{abstract}
These are notes from the mini-course given by W.~Schmid in June 2003
at the Brussels PQR2003 Euroschool.
Both authors are very thankful to Simone Gutt for organizing the
conference and her hospitality.
\end{abstract}

\tableofcontents
\newpage

\begin{section}
{Reductive Lie Groups: Definitions and Basic Properties}
\end{section}

The results stated in this section are fairly standard. Proofs and further details can be found in \cite{Kn}, for instance.

\begin{subsection}
{Basic Definitions and Examples}
\end{subsection}

In these notes ``Lie algebra" means finite dimensional Lie algebra over $\BB R$ or $\BB C$. These arise as Lie algebras of, respectively, Lie groups and complex Lie groups. We begin by recalling some basic definitions:

\begin{df}\label{defreductive}
A Lie algebra $\g g$ is {\em simple} if it has no proper ideals
and $\dim \g g >1$. A Lie algebra $\g g$ is {\em semisimple} if it can be written as a direct sum of simple ideals $\g g_i$,
$$
\g g = {\oplus}_{1\leq i\leq N}\ \g g_i\,.
$$
One calls a Lie algebra $\g g$ {\em reductive} if it can be written as a direct sum of ideals
$$
\g g = \g s \oplus \g z\,,
$$
with $\g s$ semisimple and $\g z=$ center of $\g g$. A Lie group is {\em simple}, respectively {\em semisimple}, if it has finitely many connected components and if its Lie algebra is simple, respectively semisimple. A Lie subgroup $G\subset GL(n, \BB R)$ or $G \subset GL(n, \BB C)$ is said to be {\em reductive} if it has finitely many connected components, its Lie algebra is reductive and $Z_{G^0}=$\ center of the identity component $G^0$ of $G$ consists of semisimple linear transformations~-- equivalently, if $Z_{G^0}$ is conjugate, possibly after extension of scalars from $\BB R$ to $\BB C$, to a subgroup of the diagonal subgroup in the ambient $GL(n, \BB R)$ or $GL(n, \BB C)$.
\end{df}

\begin{rem}\label{reductivegroups}
Our definition of a reductive Lie algebra is not the one most commonly used, but is equivalent to it. The semisimple ideal $\g s$ is uniquely determined by $\g g$ since $\g s = [\g g, \g g]$.

In the definition of a simple Lie algebra we require $\dim \g g >1$  because we want to exclude the one-dimensional abelian Lie algebra, which is reductive.

We shall talk about reductive Lie group only in the context of linear groups, i.e., for Lie subgroups of $GL(n, \BB R)$ or $GL(n, \BB C)$. Note that the Lie groups
$$
\left\{ \left. \begin{pmatrix}e^t & 0\\ 0 & e^{-t} \end{pmatrix}\  \right| \ t\in \BB R\ \right\}\ ,\qquad \left\{ \left. \begin{pmatrix}1 & x\\ 0 & 1 \end{pmatrix}\  \right| \ x\in \BB R\ \right\}
$$
are both isomorphic to the additive group of real numbers, and hence to each other, but only the first is reductive in the sense of our definition.
\end{rem}

\begin{ex}  \label{linear}
The Lie groups $SL(n, \BB R)$, $SL(n, \BB C)$ are simple; the Lie groups $GL(n, \BB R)$, $GL(n, \BB C)$ are reductive. Any compact real Lie group is a linear group, as can be deduced from the Peter-Weyl Theorem, and is moreover reductive (see for example Proposition 1.59, Theorem 4.20 and its Corollary 4.22 in \cite{Kn}).
\end{ex}

\begin{rem}
We will show in Example \ref{counter} that the universal covering  $\widetilde{SL(n, \BB R)}$ of $SL(n,\BB R)$, $n \ge 2$, is not a linear group.
\end{rem}

One studies reductive Lie groups because these are the groups that
naturally arise in geometry, physics and number theory. The notion of semisimple Lie group is a slight generalization of the notion
of simple Lie group, and the notion of reductive Lie group in turn generalizes the notion of semisimple Lie group.

\begin{subsection}
{Maximal Compact Subgroups and the Cartan Decomposition}
\label{cartandecomposition}
\end{subsection}

From now on we shall use the symbol $G_{\BB R}$ to denote a {\em linear reductive} Lie group; we reserve the symbol $G$ for the complexification of $G_{\BB R}$ -- see section \ref{lineargroups}. As a standing assumption, we suppose that $Z_{G_{\BB R}^0}=$\ center of the identity component $G_{\BB R}^0$ can be expressed as a direct product
\begin{equation}\label{defreductive2}
\begin{aligned}
Z_{G_{\BB R}^0} \ =\ C\cdot A\,, \ \text{with $C$ compact, $A\ \cong\ (\BB R^k,+)$ for some $k\geq 0$, and}
\\
\text{$\g a_{\BB R}\ $}=\ \text{Lie algebra of $A\,\ $ has a basis $\{\xi_1,\dots,\xi_k\}$ such that each $\xi_j$}
\\
\qquad\qquad\qquad\qquad\text{is diagonalizable, with rational eigenvalues.}
\end{aligned}
\end{equation}
The $\xi_j$ are regarded as matrices via the embedding $G_{\BB R}\subset GL(n, \BB R)$ or $G_{\BB R} \subset GL(n, \BB C)$ that exhibits $G_{\BB R}$ as linear group. In the present section, real eigenvalues would suffice, but the rationality of the eigenvalues will become important when we complexify $G_{\BB R}$.
\medskip

\begin{ex}  \label{defreductive2ex}
The groups $GL(n, \BB R)$, $GL(n, \BB C)$ satisfy the condition (\ref{defreductive2}): for $G_{\BB R}=GL(n, \BB R)$, $Z_{G_{\BB R}^0}=\{\,c \cdot 1_{n\times n}\,\mid\,c\in \BB R,\,c^n>0\}$ and $C=\{\pm 1\cdot 1_{n\times n}\}$ or $C=\{1_{n\times n}\}$, depending on whether $n$ is even or odd. For $G_{\BB R}=GL(n, \BB C)$, $Z_{G_{\BB R}^0}$ is the group of non-zero complex multiples of the identity matrix and $C=\{\,c \cdot 1_{n\times n}\,\mid\,|c|=1\}$. In both cases, $\g a_{\BB R}$ consists of all real multiples of the identity matrix, and $A=\exp(\g a_{\BB R})\cong (\g a_{\BB R},+)\cong (\BB R,+)$.
\end{ex}

We let $K_{\BB R} \subset G_{\BB R}$ denote a {\em maximal compact subgroup}. Maximal compact subgroups exist for dimension reasons. The following general facts can be found in \cite{He}, for example.
\bigskip

\begin{prop}
Under the stated hypotheses on $G_{\BB R}$,
\begin{itemize}

\item[\rm{a)}]
any compact subgroup of $G_{\BB R}$ is contained in some maximal compact subgroup $K_{\BB R}$, and $\dim_{\BB R} K_{\BB R} \ge 1$ unless $G_{\BB R}$ is abelian;

\item[\rm{b)}]
any two maximal compact subgroups of $G_{\BB R}$ are conjugate by an
element of $G_{\BB R}$;

\item[\rm{c)}]
the inclusion $K_{\BB R}\hookrightarrow G_{\BB R}$ induces an isomorphism of component groups $K_{\BB R}/ K_{\BB R}^0 \simeq G_{\BB R}/ G_{\BB R}^0$;

\item[\rm{d)}]
if $G_{\BB R}$ is semisimple, the normalizer of $K_{\BB R}$ in $G_{\BB R}$ coincides with $K_{\BB R}$.
\end{itemize}
\end{prop}
\medskip

Since the maximal compact subgroups are all conjugate,
the choice of any one of them is non-essential. At various points, we shall choose a maximal compact subgroup; the particular choice will not matter.

Let $\g g_{\BB R}$ and $\g k_{\BB R}$ denote the Lie algebras of $G_{\BB R}$ and $K_{\BB R}$, respectively. Then $K_{\BB R}$ acts on $\g g_{\BB R}$ via the restriction of the adjoint representation $Ad$. We recall the {\em Cartan decomposition} of $\g g_{\BB R}$:

\begin{prop}
There exists a unique $K_{\BB R}$-invariant linear complement $\g p_{\BB R}$ of $\g k_{\BB R}$ in $\g g_{\BB R}$ (so $\g g_{\BB R} = \g k_{\BB R} \oplus \g p_{\BB R}$ as direct sum of vector spaces), with the following properties:

\noindent {\rm a)}\ The linear map $\theta: \g g_{\BB R} \to \g g_{\BB R}$, defined by
\begin{equation*}
\theta \xi =
\begin{cases}
\xi, & \text{if $\xi \in \g k_{\BB R}$};  \\
-\xi, & \text{if $\xi \in \g p_{\BB R}$},
\end{cases}
\end{equation*}
is an {\em involutive automorphism} of $\g g_{\BB R}$; equivalently, $[\g p_{\BB R}, \g p_{\BB R}] \subset \g k_{\BB R}$
and $[\g k_{\BB R}, \g p_{\BB R}] \subset \g p_{\BB R}$.

\noindent {\rm b)}\ Every $\xi\in \g p_{\BB R}\subset \g g \g l(n,\BB C)$ is diagonalizable, with real eigenvalues; here $\g p_{\BB R} \subset \g g_{\BB R} \subset \g g \g l(n,\BB C)$ refers to the inclusion induced by the inclusion of Lie groups $G_{\BB R} \subset GL(n, \BB R) \subset GL(n, \BB C)$ or $G_{\BB R} \subset GL(n, \BB C)$.
\end{prop}

\begin{rem}
Analogously to b), every $\xi\in \g k_{\BB R}\subset \g g\g l(n,\BB C)$ is diagonalizable, with purely imaginary eigenvalues. Indeed, $\{t\mapsto \exp t\xi\}$ is a one-parameter subgroup of $K_{\BB R}$, and must therefore have bounded matrix entries; that is possible only when $\xi$ is diagonalizable over $\BB C$, with purely imaginary eigenvalues.
\end{rem}

One calls $\theta: \g g_{\BB R} \to \g g_{\BB R}$ the {\em Cartan involution} of $\g g_{\BB R}$. It lifts to an involutive automorphism of $G_{\BB R}$, which we also denote by $\theta$.

\begin{ex}
The group $G_{\BB R} = SL(n, \BB R)$ contains $K_{\BB R} = SO(n)$ as maximal compact subgroup. In this situation
\begin{align*}
\g g_{\BB R} &= \{ \xi \in \End(\BB R^n) \, \mid \, \tr \xi =0 \},  \\
\g k_{\BB R} &=
\{ \xi \in \End(\BB R^n) \, \mid \, { }^t \xi = -\xi,\: \tr \xi =0 \},  \\
\g p_{\BB R} &= \{ \xi \in \End(\BB R^n) \, \mid \, { }^t \xi = \xi,\: \tr \xi =0 \}.
\end{align*}
On the Lie algebra level, $\theta \xi = - { }^t \xi$, and on the group level, $\theta g = { }^t g^{-1}$. The group $K_{\BB R}$ can be characterized as the fixed point set of $\theta$, i.e., $K_{\BB R} = \{ g \in G_{\BB R} ;\: \theta g =g \}$.
\end{ex}

The Cartan decomposition of $\g g_{\BB R}$ has a counterpart on the group level, the so-called {\em global Cartan decomposition}:

\begin{prop} 
The map $K_{\BB R} \times \g p_{\BB R} \to G_{\BB R}$, defined by $(k, \xi) \mapsto k \cdot \exp \xi$, determines a diffeomorphism of manifolds.
\end{prop}

In the setting of the above example the proposition is essentially equivalent to the well-known assertion that any invertible real square matrix can be expressed uniquely as the product of an orthogonal matrix and a positive definite symmetric matrix.

\begin{rem}  \label{globalCartanrem}
As a consequence of the proposition, $K_{\BB R} \hookrightarrow G_{\BB R}$ is a strong deformation retract. In particular this inclusion induces isomorphisms of homology and homotopy groups.
\end{rem}

\begin{subsection}
{Complexifications of Linear Groups}
\label{lineargroups}
\end{subsection}

We continue with the hypotheses stated in the beginning of the last subsection. Like any linear Lie group, $G_{\BB R}$ has a {\em complexification} -- a complex Lie group $G$, with Lie algebra
\begin{equation}
\g g \ =_{\text{def}}\ \BB C {\otimes}_{\BB R}\, \g g_{\BB R}\,,
\end{equation}
containing $G_{\BB R}$ as a Lie subgroup, such that\newline
\indent a)\ \ the inclusion $G_{\BB R} \hookrightarrow G$ induces $\g g_{\BB R} \hookrightarrow \g g$, $\xi \mapsto 1\otimes \xi$, and\newline
\indent b)\ \ $G_{\BB R}$ meets every connected component of $G$.\newline
\noindent To construct a complexification, one regards $G_{\BB R}$ as subgroup\footnote{If $G_{\BB R}$ is presented as a linear group $G_{\BB R}\subset GL(m,\BB C)$, one uses the usual inclusion $GL(m,\BB C)\hookrightarrow GL(2m,\BB R)$ to exhibit $G_{\BB R}$ as subgroup of $GL(n,\BB R)$, with $n=2m$.} of $GL(n,\BB R)$, so that $\g g_{\BB R} \subset \g {gl}(n,\BB R)$. That makes $\g g$ a Lie subalgebra of $\g {gl}(n,\BB C) = \BB C \otimes_{\BB R}\, \g {gl}(n,\BB R)$. Then $G^0$, the connected
Lie subgroup of $GL(n,\BB C)$ with Lie algebra $\g g$, satisfies the condition a). By construction $G_{\BB R}^0\subset G^0$, and $G_{\BB R}$ normalizes $G^0$, hence $G = G_{\BB R}\cdot G^0$ is a complex Lie group with Lie algebra $\g g$, which contains $G_{\BB R}$ and satisfies both a) and b). When $G$ is a complexification of $G_{\BB R}$, one calls $G_{\BB R}$ a {\em real form} of $G$. We do not exclude the case of a Lie group $G_{\BB R}$ which happens to be a complex Lie group; in the case of $G_{\BB R}=GL(n,\BB C)$, for example, $G\cong GL(n,\BB C)\times GL(n,\BB C)\subset GL(2n,\BB C)$.

In general, the complexification of a linear Lie group depends on its realization as a linear group. In our situation, the complexification $G$ inherits the property (\ref{defreductive2}) from $G_{\BB R}$. It implies that $G$ is determined by $G_{\BB R}$ up to isomorphism, but the {\em embedding} does depend on the realization as real group, unless $G_{\BB R}^0$ has a compact center.

One can complexify the Cartan decomposition: let $\g g = \BB C \otimes_{\BB R} \g g_{\BB R}$ as before, $\g k = \BB C \otimes_{\BB R} \g k_{\BB R}$, and $\g p = \BB C \otimes_{\BB R} \g p_{\BB R}$; then $\g g = \g k \oplus \g p$ as vector spaces. The complexification $G$ of $G_{\BB R}$ naturally contains $K=$ complexification of $K_{\BB R}$, as complex Lie subgroup. Observe that $K$ cannot be compact unless $K_{\BB R}=\{e\}$, which  does not happen unless $G_{\BB R}$ is abelian; indeed, any non-zero $\xi\in\g k_{\BB R}$ is diagonalizable over $\BB C$, with purely imaginary eigenvalues, not all zero, so the complex one-parameter subgroup $\{z\mapsto \exp(z\xi)\}$ of $K$ is unbounded. By construction, the Lie algebras $\g g_{\BB R}$, $\g k_{\BB R}$, $\g g$, $\g k$ and the corresponding groups satisfy the following containments:
\begin{equation}
\begin{matrix}
\g g_{\BB R}  & \subset & \g g & \qquad\qquad & G_{\BB R}  & \subset & G \\
\cup & \qquad & \cup & \qquad\qquad & \cup & \qquad & \cup \\
\g k_{\BB R} & \subset & \g k & \qquad\qquad & K_{\BB R} & \subset & K
\end{matrix}
\end{equation}
Since $[\g p_{\BB R}, \g p_{\BB R}] \subset \g k_{\BB R}$
and $[\g k_{\BB R}, \g p_{\BB R}] \subset \g p_{\BB R}$,
\begin{equation}
\g u_{\BB R} \ =_{\text{def}}\ \g k_{\BB R} \oplus i \g p_{\BB R}\ \ \ \text{is a real Lie subalgebra of $\g g$.}
\end{equation}
We let $U_{\BB R}^0$ denote the connected Lie subgroup of $G$ with Lie algebra $\g u_{\BB R}$. If $G_{\BB R}$ is semisimple, one knows that $U_{\BB R}^0$ is compact \cite{He}; as a consequence of our hypotheses, $U_{\BB R}^0$ is compact even in the reductive case. Thus $U_{\BB R}^0$ lies in a maximal compact subgroup of $G$, which we denote by $U_{\BB R}$. Not only is $U_{\BB R}\subset G$ a maximal compact subgroup, but also
\medskip
\begin{equation}
\begin{aligned}
\rm{a)} \ \ &U_{\BB R} \, \ \text{is a real form of}\,\ G;\qquad\qquad\qquad\qquad\qquad\qquad\qquad\qquad\qquad
\\
\rm{b)}\ \ &K_{\BB R} \ = \ U_{\BB R} \cap G_{\BB R}.
\end{aligned}
\end{equation}
Both assertions are well known in the semisimple case, to which the general case can be reduced.

The process of complexification establishes a bijection, up to isomorphism, between compact Lie groups and linear, reductive, complex Lie groups\footnote{Any connected complex Lie group with a reductive Lie algebra can be realized as a linear group; we mention the hypothesis of linearity to signal that we want the linear realization to be reductive in the sense of our definition.} satisfying (\ref{defreductive2}); in the opposite direction, the correspondence is given by the passage to a maximal compact subgroup, which is then a {\em compact real form} of the complex group. The groups $U_{\BB R}$ and $G$ are related in this fashion: the former is a compact real form of the latter. In view of Remark \ref{globalCartanrem}, $U_{\BB R} \hookrightarrow G$ is a strong deformation retract, which induces isomorphisms of homology and homotopy.
\medskip

\begin{ex}
The pair $G_{\BB R}=SL(n, \BB R)$, $K_{\BB R} = SO(n)$, has complexifications $G = SL(n, \BB C)$, $K = SO(n, \BB C)$; the corresponding compact real form of $\,G = SL(n, \BB C)$ is $U_{\BB R} = SU(n)$.
\end{ex}

Since $\g g = \BB C\otimes_{\BB R}\,\g u_{\BB R}$, these two Lie algebras have the same representations over $\BB C$ -- representations can be restricted from $\g g$ to $\g u_{\BB R}$, and in the opposite direction, can be extended complex linearly. On the global level, these operations induce a canonical bijection
\begin{equation}\label{extensionofrep}
\left\{ \begin{matrix} \text{holomorphic finite dimensional} \\ \text{representations of $G$}\end{matrix}\right\}
\ \ \simeq \ \
\left\{ \begin{matrix} \text{finite dimensional continuous} \\ \text{complex representations of $U_{\BB R}$}\end{matrix}\right\}\,,
\end{equation}
a bijection because $U_{\BB R}\hookrightarrow G$ induces an isomorphism of the component group and the fundamental group. Of course we are also using the well known fact that continuous finite dimensional representations of Lie groups are necessarily real analytic, and are determined on the identity component by the corresponding infinitesimal representations of the Lie algebra.

We had mentioned earlier that the universal covering group of $G=SL(n, \BB R)$, $n \geq 2$, is not a linear group. We can now sketch the argument:

\begin{ex}  \label{counter}
Let $\widetilde{SL(n,\BB R)}$ be the universal covering group of $SL(n,\BB R)$, $n \ge 2$. Since
$$
\pi_1(SL(n,\BB R)) = \pi_1(SO(n)) =
\begin{cases}
\BB Z,        & \text{if $n=2$;}  \\
\BB Z/2\BB Z, & \text{if $n \ge 3$,}
\end{cases}
$$
the universal covering $\widetilde{SL(n,\BB R)} \rightarrow SL(n,\BB R)$ is a principal $\BB Z$-bundle when $n=2$ and a principal $\BB Z/2\BB Z$-bundle when $n \ge 3$. If $\widetilde{SL(n,\BB R)}$ were linear, its complexification would have to be a covering group of
$SL(n,\BB C)=$ complexification of $SL(n, \BB R)$, of infinite order when $n=2$ and of order (at least) two when $n \ge 3$. But $SU(n)$, $n \ge 2$, is simply connected, as can be shown by induction on $n$. But then $SL(n,\BB C)=$ complexification of $SU(n)$ is also simply connected, and therefore cannot have a non-trivial covering. We conclude that $\widetilde{SL(n,\BB R)}$ is not a linear group.
\end{ex}

\begin{section}
{Compact Lie Groups}
\end{section}

In this section we consider the well understood case of a connected, compact Lie group. As was remarked in Example \ref{linear}, any such group is automatically linear and reductive. In the setting of chapter 1, the groups $K_{\BB R}$ and $U_{\BB R}$ are compact, but not necessarily connected. In any case, knowing the representations of the identity component explicitly gives considerable information about the representations of a non-connected compact group -- modulo knowledge of the representations of the component group, of course. In section 4, we shall suppose that the group $G_{\BB R}$ has a connected complexification $G$; in that case, $U_{\BB R}$ will indeed be connected. Let us state the hypotheses of the current section explicitly:
\begin{equation}
U_{\BB R} \ \ \text{is a connected compact Lie group}.
\end{equation}
As in the previous section $\g g=\BB C \otimes_{\BB R}\,u_{\BB R}$ denotes the complexified lie algebra of $U_{\BB R}$.

\bigskip

\begin{subsection}
{Maximal Tori, the Unit Lattice, and the Weight Lattice}
\end{subsection}

A common strategy in mathematics is to study properties of a
new class of objects by looking for sub-object whose properties are already known, but which are ``large enough" to convey some useful information about the objects to be studied. In representation theory this means studying representations of compact groups by restricting them to maximal tori, while representations of noncompact linear reductive groups are studied by restricting them to maximal compact subgroups.

With $U_{\BB R}$ connected and compact, as we are assuming, let $T_{\BB R} \subset U_{\BB R}$ denote a {\em maximal torus}. It is not difficult to see that maximal tori exist and are nontrivial -- i.e., $\dim T_{\BB R} > 0$ -- unless $U_{\BB R}=\{e\}$. Moreover,
\medskip 

\begin{prop}\label{conjugacyoftori}\qquad
\begin{itemize}
\item[\rm{a)}]
Any two maximal tori in $U_{\BB R}$ are conjugate by an element of $U_{\BB R}$;
\item[\rm{b)}]
any $g \in U_{\BB R}$ is conjugate to some $t \in T_{\BB R}$;
\item[\rm{c)}]
$T_{\BB R}$ coincides with its own centralizer in $U_{\BB R}$, i.e., $T_{\BB R} = Z_{U_{\BB R}}(T_{\BB R})$;
\item[\rm{d)}]
$T_{\BB R}$ is the identity component of its own normalizer in $U_{\BB R}$, i.e., $T_{\BB R} = N_{U_{\BB R}}(T_{\BB R})^0$.
\end{itemize}
\end{prop}

We fix a particular maximal torus $T_{\BB R}$, with Lie algebra $\g t_{\BB R}$ and complexified Lie algebra $\g t$. In view of the Proposition, the particular choice will not matter. Since $T_{\BB R}$ is abelian and connected, the exponential $\exp : \g t_{\BB R} \to T_{\BB R}$ is a surjective homomorphism, relative to the additive structure of $\g t_{\BB R}$.

\begin{rem}
The exponential mapping of a general connected Lie group need not be surjective. For instance, $g = \begin{pmatrix} -1 & 1 \\ 0 & -1 \end{pmatrix} \in SL(2,\BB R)$ cannot lie in the image of the exponential map. Indeed, if this element of $g$ could be expressed as $g=\exp \xi$, for some $\xi \in \g{sl} (2,\BB R)$, then $\xi$ is not diagonalizable, even over $\BB C$, since $g$ is not diagonalizable over $\BB C$. That forces $\xi$ to have two equal eigenvalues, necessarily eigenvalues zero because $\tr \xi =0$. Contradiction: $g=\exp(\xi)$ does not have eigenvalues $1$.
\end{rem}

The exponential map $\exp : \g t_{\BB R} \to T_{\BB R}$ is not only a surjective homomorphism, but also locally bijective, hence a covering homomorphism,
\begin{equation}
\exp : \g t_{\BB R} /L \ \tilde{\longrightarrow} \ T_{\BB R}\qquad (\, L = \{ \xi \in t_{\BB R} \,\mid\, \exp \xi = e \}\,)\,.
\end{equation}
That makes $L\subset t_{\BB R}$ a discrete, cocompact subgroup. In other words, $L$ is a lattice, the so-called {\em unit lattice}. Let $\widehat T_{\BB R}$ denote the group of characters, i.e. the group of homomorphisms from $T_{\BB R}$ to the unit circle $\,S^1 = \{ z \in \BB C \mid |z|=1 \}$. Then
\begin{equation}
\Lambda \ =_{\text{def}} \ \{\, \lambda \in i \g t_{\BB R}^* \,\mid\, \langle \lambda, L \rangle \subset 2\pi i \BB Z \,\} \ \ \overset{\sim}{\longrightarrow} \ \ \widehat T_{\BB R}\,, \qquad \Lambda \, \ni \, \lambda\, \longmapsto \, e^\lambda \, \in \, \widehat T_{\BB R}\,,
\end{equation}
with $e^{\lambda} : T_{\BB R} \to \{ z \in \BB C \,\mid \,|z|=1 \}$
defined by $e^{\lambda}(\exp \xi) = e^{\langle \lambda, \xi \rangle}$, for any $\xi \in \g t_{\BB R}$. One calls $\Lambda \subset i \g t_{\BB R}^*$ the {\em weight lattice}; except for the factor $2\pi i$ in its definition, it is the lattice dual of the unit lattice $L\subset\g t_{\BB R}$.

\begin{subsection}
{Weights, Roots, and the Weyl Group}
\end{subsection}

Let $\pi$ be a representation of $U_{\BB R}$ on a finite-dimensional complex vector space $V$ -- in other words, a continuous homomorphism $\pi: U_{\BB R} \to GL(V)$. Since $T_{\BB R}$ is compact and the field $\BB C$ algebraically closed, the action of any $t\in T_{\BB R}$ can be diagonalized, and since $T_{\BB R}$ is abelian, the action of all the $t \in T_{\BB R}$ can be diagonalized simultaneously. Thus we obtain the {\em weight space decomposition}
\begin{equation}
V \ =\  {\oplus}_{\lambda \in \Lambda} \,V^{\lambda}\,,
\end{equation}
where
\begin{equation}
V^{\lambda} \ = \ \{\, v \in V \,\mid\, \pi(t)v = e^{\lambda}(t)v \,\ \forall t \in T_{\BB R} \,\} \ = \ \{\, v \in V \,\mid\, \pi(\xi)v = \langle \lambda, \xi \rangle v \,\ \forall \xi \in \g t_{\BB R}\, \}\,.
\end{equation}
If $V^{\lambda} \ne \{0\}$, one calls $\lambda$, $V^{\lambda}$, and
$\dim V^{\lambda}$ respectively a {\em weight} of $\pi$, the {\em weight space} corres\-ponding to $\lambda$, and the {\em multiplicity} of the weight $\lambda$.

In the case of the adjoint representation of $U_{\BB R}$ on the complexified Lie algebra $\g g=\BB C \otimes_{\BB R}\g u_{\BB R}$, one singles out the weight zero: $\g g = \g g^0 \oplus \left( {\oplus}_{\lambda \ne 0}\, \g g^{\lambda} \right)$. Evidently $\g g^0 \supset \g t$, since $T_{\BB R}$ acts trivially on $\g t$. In fact, $\g g^0 = \g t$, for one could otherwise show that $T_{\BB R}$ is not a maximal torus. Nonzero weights of the adjoint representation are called {\em roots}, hence
\begin{equation}
\g g\  = \ \g t\ \oplus \ \left( {\oplus}_{\alpha \in \Phi}\, \g g^{\alpha} \right),\ \ \ \text{with}\ \ \Phi \,=\, \Phi(\g g,\g t)\, = \, \text{set of all roots}\,.
\end{equation}
This is the {\em root space decomposition} of $\g g$ relative to the action of $T_{\BB R}$, or equivalently, relative to the action of $\g t$. One refers to $\Phi$ as the {\em root system}. Very importantly,
\begin{equation}
\alpha \in \Phi \ \Longrightarrow \ \dim \g g^\alpha \ = \ 1\,,
\end{equation}
i.e., all roots have multiplicity one. Since $\Phi \subset\Lambda - \{0\}\subset i \g t_{\BB R}^* \subset \g t^*$, roots take purely imaginary values on the real Lie algebra $\g t_{\BB R}$, which implies
\begin{equation}
\overline{\g g^{\alpha}} \ =\ \g g^{\overline \alpha}\ = \ \g g^{-\alpha}\,,\ \ \text{and hence}\ \ \Phi \ = \ - \Phi\,;
\end{equation}
here $\overline{\g g^{\alpha}}$ denotes the complex conjugate of $\g g^{\alpha}$ with respect to the real form $\g u_{\BB R}\subset \g g$.

Let $\pi : \g u_{\BB R} \to \g g \g l(V)$ denote the infinitesimal representation induced by the global representation $\pi$ considered at the beginning of this subsection, and $\pi : \g g \to \g g \g l(V)$ its complex extension. The latter may be interpreted as a linear map $\g g\otimes_{\BB C} V \to V$, which is $U_{\BB R}$-invariant -- therefore also $T_{\BB R}$-invariant -- when $U_{\BB R}$ is made to act on $\g g$ via $\operatorname{Ad}$. Hence, for every $\alpha\in\Phi$ and every weight $\lambda$ of $\pi$,
\begin{equation}
\pi(\g g^\alpha)\,V^\lambda \ \ \subset \ V^{\lambda + \alpha}\,; \ \ \ \text{in particular}\ \ \pi(\g g^\alpha)\,V^\lambda \, = \, 0\ \ \ \text{if $\lambda + \alpha$ is not a weight}.
\end{equation}
Applied to the adjoint representation, this means
\begin{equation}
[ \g g^\alpha , \g g^\beta ] \ \subset \ 
\begin{cases}
\ \g g^{\alpha+\beta}\ \ &\text{if}\,\ \alpha + \beta \in \Phi
\\
\ \g t\ \ &\text{if}\,\ \alpha + \beta = 0
\\

\ 0\ \ &\text{if}\,\ \alpha + \beta \notin \Phi \cup \{0\}\,,
\end{cases}
\end{equation}
for all roots $\alpha,\beta\in\Phi$.

An element $\xi\in\g t$ is said to be {\em singular} if $\langle \alpha,\xi\rangle = 0$ for some root $\alpha$, and otherwise {\em regular}.
The set
\begin{equation}
i\g t'_{\BB R}\ =_{\text{def}} \ \{\,\g it_{\BB R}\ \mid \ \text{$\xi$ is regular}\,\}
\end{equation}
breaks up into a finite, disjoint union of open, convex cones, the so-called {\em Weyl chambers}. If $C \subset i \g t_{\BB R}$ is a particular Weyl chamber and $\xi$ an element of $C$, the subset
\begin{equation}
\Phi^+ \ = \ \{ \alpha \in \Phi \,\mid\, \langle \alpha, \xi \rangle >0 \}\ \subset \ \Phi
\end{equation}
depends only on $C$, not on the choice of $\xi \in C$; by definition, $\Phi^+$ is a {\em system of positive roots}. Essentially by construction,
\begin{equation}
\begin{aligned}
&\ \rm{a)}\ \ \ \ \Phi\ = \ \Phi^+ \cup (-\Phi^+)\ \ \ \text{(disjoint union);}
\\
&\ \rm{b)}\ \ \ \ \alpha, \beta \in \Phi^+\,,\ \ \alpha + \beta \in \Phi \ \ \Longrightarrow \ \ \alpha + \beta \in \Phi^+.
\end{aligned}
\end{equation}
The Weyl chamber $C$ can be recovered from the system of positive roots $\Phi^+$,
\begin{equation}
C \ =  \ \{\, \xi \in i \g t_{\BB R} \, \mid\, \langle \alpha, \xi \rangle >0 \,\}\,.
\end{equation}
In fact, $C \,\longleftrightarrow\, \Phi^+$ establishes a bijection between Weyl chambers and positive root systems.

Via the adjoint action, the normalizer $N_{U_{\BB R}}(T_{\BB R})$ acts on $\g t_{\BB R}$, on $i\g t_{\BB R}^*$, on $\Lambda$, on $\Phi$, and on the set of Weyl chambers. Since $T_{\BB R} = Z_{U_{\BB R}}(T_{\BB R}) =  N_{U_{\BB R}}(T_{\BB R})^0$ by Proposition \ref{conjugacyoftori}, the {\em Weyl group}
\begin{equation}
W\ = \ W(U_{\BB R}, T_{\BB R})\ =_{def}\ N_{U_{\BB R}}(T_{\BB R}) /T_{\BB R}
\end{equation}
is a finite group, which acts on $\g t_{\BB R}$, on $i\g t_{\BB R}^*$, on the weight lattice $\Lambda$, and on $\Phi$. This action permutes the Weyl chambers, hence also the positive root systems.
\medskip

\begin{prop} \label{simplytransitively}
The Weyl group $W(U_{\BB R}, T_{\BB R})$ acts faithfully on $\g t_{\BB R}$, and acts simply transitively on the set of Weyl chambers $\{C \}$, as well as on the set of positive root systems $\{\Phi^+\}$.
\end{prop}

In particular, if some $n\in N_{U_{\BB R}}(T_{\BB R})$ fixes a root system $\Phi^+$, then $n\in T_{\BB R}$. In the following, we shall fix a positive root system $\Phi^+$. The particular choice will not matter since they are all conjugate to each other.  The Weyl chamber $C$ that corresponds to $\Phi^+$ is called the {\em dominant Weyl chamber}.

\begin{subsection}
{The Theorem of the Highest Weight}
\end{subsection}

We consider finite dimensional representations of $U_{\BB R}$ on complex vector spaces, as in the previous subsection. Recall that a finite dimensional representation $(\pi,V)$ is {\em irreducible} if the representation space $V$ contains no $U_{\BB R}$-invariant subspaces other than $\{0\}$ and $V$ itself; it is {\em completely reducible} if it can be written as the direct sum of irreducible subrepresentations. If $(\pi, V)$ is unitary -- i.e., if $V$ comes equipped with $U_{\BB R}$-invariant inner product -- the orthogonal complement of an invariant subspace is again invariant. One can therefore successively split off one minimal invariant subspace at a time. Since minimal invariant subspaces are irreducible, this shows that finite dimensional, unitary representations are completely reducible. Any representation $(\pi,V)$ of the compact group $U_{\BB R}$ can be made unitary: one puts an arbitrary inner product on the space $V$, and then uses Haar measure to average the $g$-translates of this inner product for all $g\in U_{\BB R}$; the averaged inner product is $U_{\BB R}$-invariant. This implies the well known fact that
\begin{equation}
\text{finite dimensional representations of compact groups are completely reducible}.
\end{equation}    
In particular, to understand the finite dimensional representations of $U_{\BB R}$, it suffices to understand the finite dimensional, irreducible representations. 

The preceding discussion applies to any compact Hausdorff group. We now return to the case of a connected, compact Lie group $U_{\BB R}$. Our definitions and statements will involve the choice of a maximal torus $T_{\BB R}$ and a positive root system $\Phi^+ \subset i \g t_{\BB R}^*$. As was remarked earlier, these are not essential choices. Since $U_{\BB R}$ is compact, there exists a negative definite, $\Ad(U_{\BB R})$-invariant, symmetric bilinear form
\begin{equation}\label{Sdef}
S\, : \, \g u_{\BB R} \times \g u_{\BB R} \ \longrightarrow \ \BB R\,,
\end{equation}
$\Ad(U_{\BB R})$-invariant in the sense that $S \bigl( \Ad\, g(\xi), \Ad\,g(\eta) \bigr) = S (\xi, \eta)$, for all $\xi,\eta \in \g u_{\BB R}\,$ and  $g \in U_{\BB R}\,$, or equivalently, on the infinitesimal level,
\begin{equation}
S \bigl( \,[\zeta, \xi]\,,\, \eta \,\bigr) + S \bigl(\, \xi\,,\, [\zeta, \eta] \,\bigr) \ = \ 0\,, \qquad \text{for all}\,\ \zeta, \xi, \eta \in \g u_{\BB R}\,.
\end{equation}
One way to construct $S$ is to take a finite-dimensional representation $(\pi,V)$ which is faithful, or at least faithful on the level of the Lie algebra, and define
\begin{equation}
S(\xi,\eta)\ =\ \tr \bigl( \pi(\xi)\pi(\eta) \bigr)\,.
\end{equation}
Indeed, $S$ is symmetric, bilinear, $\Ad(U_{\BB R})$-invariant by construction, and any nonzero $\xi\in\g u_{\BB R}$ acts on $V$ dagonalizably, with purely imaginary eigenvalues, not all zero, so $\tr \bigl( \pi(\xi)\pi(\xi) \bigr)<0$. If $U_{\BB R}$ is semisimple, one can let the adjoint representation play the role of $\pi$; in that case one calls $S$ the {\em Killing form}.

The bilinear form $S$ is far from uniquely determined by the required properties. However, its restriction to the various simple ideals in $\g u_{\BB R}$ {\em is unique}, up to scaling, as follows from Schur's lemma. This partial uniqueness suffices for our purposes.

Extending scalars from $\BB R$ to $\BB C$, we obtain an $\Ad(U_{\BB R})$-invariant, symmetric, complex bilinear form $S : \g g \times \g g \to \BB C$, which is positive definite on $i\g u_{\BB R}$. By restriction it induces a positive definite inner product
$(\,.\,,\,.\,)$ on $i \g t_{\BB R}$, and by duality also on $i \g t_{\BB R}^*$. The Weyl group $W=W(U_{\BB R},T_{\BB R})$ preserves these inner products, since they are obtained by restriction of an $Ad(U_{\BB R})$-invariant bilinear form.

\medskip

\begin{df}  \label{dominant}
An element $\lambda\in i\g t_{\BB R}^*$ is said to be {\em dominant} if $(\lambda, \alpha) \ge 0$ for all $\alpha \in \Phi^+$, and {\em regular} if $(\lambda, \alpha) \ne 0$ for all $\alpha \in \Phi$.
\end{df}

These notions apply in particular to any $\lambda$ in the weight lattice $\Lambda$. The inner product identifies $i\g t_{\BB R}^*$ with its own dual $i\g t_{\BB R}$, and via this identification, the set of all dominant, regular $\lambda\in i\g t_{\BB R}^*$ corresponds precisely to the dominant Weyl chamber $C\subset \g t_{\BB R}$, i.e., the Weyl chamber determined by $\Phi^+$. Since $W$ acts simply transitively on the set of Weyl chambers, every regular $\lambda\in i\g t_{\BB R}^*$ is $W$-conjugate to exactly one dominant, regular $\lambda'\in i\g t_{\BB R}^*$. In fact, this statement remains correct if one drops the condition of regularity; this can be seen by perturbing a singular $\lambda\in i\g t_{\BB R}^*$ slightly, so as to make it regular. The action of $W$ preserves the weight lattice, hence
\begin{equation}
\text{every $\lambda\in\Lambda$ is $W$-conjugate to a unique dominant $\lambda'\in\Lambda$}\,;
\end{equation}
in other words, $\{\lambda\in\Lambda\,\mid\,\text{$\lambda$ is dominant}\}\cong W\backslash \Lambda$.

\begin{thm} [Theorem of the Highest Weight]
For an irreducible, finite dimensional, complex representation $\pi$, the following conditions on a weight $\lambda$ of $\pi$ are equivalent:
\begin{enumerate}
\item
$\lambda+\alpha$ is not a weight, for any positive root $\alpha \in \Phi^+$;
\item
there exists a non-zero $v_0 \in V^{\lambda}$ such that
$\,\pi(\g g^{\alpha})\, v_0 =0$ for all $\alpha \in \Phi^+$;
\item
any weight of $\pi$ can be expressed as $\lambda - A$, where $A$ is a sum of positive roots (possibly empty; repetitions are allowed).
\end{enumerate}
There exists exactly one weight $\lambda$ of $\pi$ with these (equivalent) properties, the so-called {\em highest weight} of $\pi$. The highest weight is dominant, has multiplicity one (i.e. $\dim V^{\lambda} =1$), and determines the representation $\pi$ up to isomorphism. Every dominant $\lambda\in\Lambda$ arises as the highest weight of an irreducible representation $\pi$.
\end{thm}

The assertion that property {\it 2} implies property {\it 3} can be deduced from the Poincar\'e-Birkhoff-Witt theorem, and the other implications among the three properties can be established by elementary arguments. Among the remaining statements, only the existence of an irreducible finite dimensional representation with a given regular dominant weight requires some effort. One can prove this existence statement analytically, via the Weyl character formula, algebraically, by realizing the representation in question as a quotient of a Verma modules, or geometrically, as will be sketched in subsection \ref{The Borel-Weil-Bott Theorem}.

In effect, the theorem parameterizes the isomorphism classes of irreducible finite dimensional representations over $\BB C$ in terms of their highest weights,
\begin{equation}\label{enumerationofreps}
\Bigl\{ \begin{matrix}
\text{irreducible finite-dimensional complex} \\
\text{representations of $U_{\BB R}$, up to isomorphism}
\end{matrix} \Bigr\}
\longleftrightarrow
\{ \lambda \in \Lambda \mid \text{$\lambda$ is dominant} \} \longleftrightarrow W \backslash \Lambda\,.
\end{equation}
Beyond the enumeration of irreducible finite dimensional representations, the theorem also provides structural information about these representations. In fact, most of the general structural properties that are used in applications are consequences of the Theorem of the Highest Weight.

\begin{subsection}
{Borel Subalgebras and the Flag Variety}
\end{subsection}

We now also consider the complexification $G$ of the connected compact Lie group $U_{\BB R}$. In principle, its Lie algebra $\g g = \BB C \otimes_{\BB R}\,\g u_{\BB R}$ can be any reductive Lie algebra over $\BB C$, and $G$ any connected, complex, reductive Lie group.

\bigskip 

\begin{df}  \label{Borelsubalgebra}
A {\em Borel subalgebra} of $\g g$ is a maximal solvable subalgebra. A {\em Borel subgroup} of $G$ is a connected complex Lie subgroup of
$G$ whose Lie algebra is a Borel subalgebra of $\g g$.
\end{df}

\begin{prop}
Any two Borel subgroups of $G$, respectively Borel subalgebras of
$\g g$, are conjugate under the action of $G$. Any Borel subgroup $B \subset G$ coincides with its own normalizer, i.e., $N_G(B)=B$. In particular, Borel subgroups are closed.
\end{prop}

\begin{rem}
The property $N_G(B)=B$ implies that {\em any} complex subgroup of $G$ whose Lie algebra is a Borel subalgebra is automatically connected, and hence is a Borel subgroup. In other words, the requirement of connectedness in Definition \ref{Borelsubalgebra} can be dropped without changing the notion of a Borel subgroup.
\end{rem}

The linear subspaces of $\g g$ of a fixed dimension constitute a smooth projective variety, a so-called Grassmann variety. Being a subalgebra, or more specifically a solvable subalgebra, amounts to an algebraic condition on an arbitrary point in this Grassmannian. Since all Borel subalgebras have the same dimension, call it $d$, a solvable subalgebra of $\g g$ is maximal solvable if and only if it has dimension $d$. We conclude that
\begin{equation}
X \ = \ \text{set of all Borel subalgebras of $\g g$}
\end{equation}
is a closed subvariety of the Grassmannian. This gives $X$ the structure of a complex projective variety. By definition, $X$ is the {\em flag variety} of $\g g$. We already know that $G$ acts transitively on $X$ via the adjoint action, which is algebraic. In particular $X$ is smooth.

The flag variety $X$ can be characterized by a universal property: it dominates all the complex projective varieties with a transitive, algebraic action of $G$ -- any other variety $Y$ with these properties is a $G$-equivariant quotient of $X$, i.e., the image of $X$ under a surjective algebraic map that relates the $G$-actions on $X$ and $Y$. Such $G$-equivariant quotients of $X$ are called {\em generalized flag varieties}

As before, let $T_{\BB R}\subset U_{\BB R}$ be a maximal torus, $\g t$ its complexified Lie algebra, and $\Phi^+$ a system of positive roots. One can show quite directly that
\begin{equation}
\g b_0\ =\ \g t\ \oplus \ \left( {\oplus}_{\alpha \in \Phi^+}\, \g g^{-\alpha} \right)
\end{equation}
is maximal solvable in $\g g$, hence a Borel subalgebra. Any other Borel subalgebra is conjugate to it under the action of $G$, and even under the action of $U_{\BB R}$, as we shall see. The corresponding Borel subgroup $B_0$ is also the normalizer of $B_0$, hence the normalizer of $\g b_0$ in $G$ -- in other words, $B_0$ is the isotropy subgroup at $\g b_0$ for the action of $G$ on $X$. That implies $X \simeq G / B_0$, since $G$ acts transitively on $X$.

\begin{lem}
$U_{\BB R} \cap B_0 = T_{\BB R}$.
\end{lem}

\pf
Complex conjugation with respect to the real form $\g u_{\BB R}\subset \g g$ maps $\g g^{-\alpha}$ to $\g g^\alpha$, hence $\g b_0 \cap \overline{\g b_0}=\g t$, hence
\begin{equation}
\g u_{\BB R} \cap \g b_0\ = \ (\g u_{\BB R} \cap \g b_0) \cap (\overline{\g u_{\BB R} \cap \g b_0})\ = \ \g u_{\BB R} \cap \g b_0 \cap \overline{\g b_0}\ = \ \g u_{\BB R} \cap \g t \ = \ \g t_{\BB R}\,.
\end{equation} 
On group level that means $(U_{\BB R} \cap B_0)^0 = T_{\BB R}$. Any $n \in U_{\BB R} \cap B_0$ therefore normalizes $T_{\BB R}$; $n$ maps $\Phi^+$ to itself, since otherwise $\operatorname{Ad}n(\g b_0)$ could not equal $\g b_0$, as it must. We had mentioned earlier that any $n \in N_{U_{\BB R}}(T_{\BB R})$ which fixes $\Phi^+$ lies in $T_{\BB R}$. The lemma follows.
\medskip

Because of the lemma, we can identify the $U_{\BB R}$-orbit through the identity coset in $G/B_0\simeq X$ with $U_{\BB R}/T_{\BB R}$. This orbit is closed because $U_{\BB R}$ is compact, and is open, as can be seen by counting dimensions. Hence $U_{\BB R}$ acts transitively on $X$, and
\begin{equation}
X \ \simeq \ G / B_0\ \simeq\ U_{\BB R} / T_{\BB R} \,.
\end{equation}
The transitivity of the $U_{\BB R}$-action now implies that each Borel subgroup of $G$ intersects $U_{\BB R}$ in a maximal torus. Arguing as in the proof of the lemma, one finds that $W = N_{U_{\BB R}}(T_{\BB R})/T_{\BB R}$ acts simply transitively on the set of Borel subgroups which contain $T_{\BB R}$. Thus each maximal torus in $U_{\BB R}$ lies in exactly $N$ Borel subgroups, $N=$ cardinality of $W$.

\begin{ex}
A {\em complete flag} in $\BB C^n$ is a nested sequence of subspaces
\begin{equation*}
0 \subset F_1 \subset F_2 \subset \dots \subset F_{n-1} \subset \BB C^n\,,\qquad \,\dim F_j =j\,.
\end{equation*}
The tautological action of $G = SL(n, \BB C)$ on $\BB C^n$ induces a transitive action on the set of all such complete flags. As a consequence of Lie's theorem, any Borel subgroup of $SL(n,\BB C)$ is the stabilizer of a complete flag. In the case of $G = SL(n, \BB C)$, then, the flag variety $X$ is the variety of all complete flags in $\BB C^n$ -- hence the name {\em flag variety}.
\end{ex}

\begin{subsection}
{The Borel-Weil-Bott Theorem}\label{The Borel-Weil-Bott Theorem}
\end{subsection}

Recall the notion of a $G$-equivariant holomorphic line bundle over a complex manifold with a holomorphic $G$-action -- in our specific case, a $G$-equivariant holomorphic line bundle over the flag variety $X$: it is a holomorphic line bundle ${\cal L} \to X$, equipped with a holomorphic action of $G$ on ${\cal L}$ by bundle maps, which lies over the $G$-action on the base $X$. The isotropy group at any point $x_0\in X$ then acts on the fibre ${\cal L}_{x_0}$ at $x_0$. In this way, one obtains a holomorphic representation $\phi: B_0 \to GL(1,\BB C) = \BB C^*$ of $B_0$, the isotropy group at the identity coset in $G/B_0 \simeq X$. One dimensional representations are customarily called characters. Since $G$ acts transitively on $X$, the passage from $G$-equivariant holomorphic line bundles -- taken modulo isomorphism, as usual -- to holomorphic characters $\phi : B_0 \to \BB C^*$ can be reversed,  
\begin{equation}
\left\{\begin{matrix}
\text{holomorphic $G$-equivariant}  \\
\text{line bundles over $X \simeq G/B_0$} \end{matrix}\right\} \ \
\simeq \ \
\left\{\begin{matrix} \text{holomorphic}  \\
\text{characters of $B_0$} \end{matrix}\right\}\,.
\end{equation}
By construction this is an isomorphism of groups, relative to operations of tensor product of $G$-equivariant holomorphic line bundles and multiplication of holomorphic characters, respectively.

Holomorphic characters $\phi : B_0 \to \BB C^*$ drop to $B_0/[B_0,B_0]$, the quotient of $B_0$ modulo its commutator. Note that $B_0$ contains $T$, the complexification of the maximal torus $T_{\BB R}$. One can show that the inclusion $T \hookrightarrow B_0$ induces an isomorphism $T \simeq B_0/[B_0,B_0]$. Thus $B_0$ and $T$ have the same group of holomorphic characters. On the other hand,
\begin{equation}
\left\{\begin{matrix} \text{holomorphic}  \\
\text{characters of $T$} \end{matrix}\right\} \ \ 
\simeq \ \ \widehat{T_{\BB R}}
\end{equation}
by restriction from $T$ to its compact real form $T_{\BB R}\,$ (\ref{extensionofrep}). Combining these isomorphisms and identifying the dual group $\widehat{T_{\BB R}}$ with the weight lattice $\Lambda$ as usual, we get a canonical isomorphism
\begin{equation}
\left\{\begin{matrix}\label{Llambda}
\text{group of holomorphic $G$-equivariant}  \\
\text{line bundles over $X$} \end{matrix}\right\}\ \ 
\simeq \ \ \Lambda\,.
\end{equation}
We write ${\cal L}_{\lambda}$ for the line bundle corresponding to $\lambda \in \Lambda$ under this isomorphism.

The action of $G$ on $X$ and ${\cal L}_{\lambda}$ determines a holomorphic, linear action on the space of global sections $H^0(X, {\cal O}({\cal L}_{\lambda}))$ and, by functorality, also on the higher cohomology groups $H^p(X, {\cal O}({\cal L}_{\lambda}))$, $p>0$. These groups are finite dimensional since $X$ is compact. The Borel-Weil-Bott theorem describes the resulting representations of the compact real form $U_{\BB R}\subset G$ and, in view of (\ref{extensionofrep}), also as holomorphic representations of $G$.
\bigskip

\begin{thm} [Borel-Weil \cite{Se}]
If $\lambda$ is a dominant weight, the representation of $U_{\BB R}$ on $H^0(X, {\cal O}({\cal L}_{\lambda}))$ is irreducible, of highest weight $\lambda$, and $H^p(X, {\cal O}({\cal L}_{\lambda}))=0$ for $p>0$. If $\lambda\in\Lambda$ fails to be dominant, $H^0(X, {\cal O}({\cal L}_{\lambda}))=0$.
\end{thm}

In particular the theorem provides a concrete, geometric realization of every finite dimensional irreducible representation of $U_{\BB R}$. The description of $H^0(X, {\cal O}({\cal L}_{\lambda}))$ can be deduced from the theorem of the highest weight, and the vanishing of the higher cohomology groups is a consequence of the Kodaira vanishing theorem.

Bott \cite{Bott} extended the Borel-Weil theorem by identifying the higher cohomology groups as representations of $U_{\BB R}$. The description involves
\begin{equation}
\rho \ =\ \frac 12\ {\sum}_{\alpha \in \Phi^+}\ \alpha \ \in \ i\g t_{\BB R}^*\,.
\end{equation}
Since $\Phi\subset\Lambda$, $2\rho$ is evidently a weight. In fact,
\begin{equation}
{\cal L}_{-2\rho}\ \ = \ \ \text{canonical bundle of $X$}
\end{equation}
as can be shown quite easily. In general, $\rho$ itself need not be a weight; if not a weight, it can be made a weight by going to a twofold covering group. Geometrically this means that the canonical line bundle of $X$ has a square root, possibly as a $G$-equivariant holomorphic line bundle, and otherwise as an equivariant holomorphic line bundle for a twofold covering of $G$. Whether or not $\rho$ is a weight,
\begin{equation}
w\rho \, - \, \rho \ \in \ \Lambda\ \ \ \text{for all}\,\ w\in W\,.
\end{equation}
Also, $\rho$ has the following important property:
\begin{equation}
\text{for}\,\ \lambda \in\Lambda\,,\ \ \ \lambda\, \ \text{is dominant}\ \ \Longleftrightarrow\ \ \lambda + \rho\, \ \text{is dominant regular}.
\end{equation}
It follows that for $\lambda\in\Lambda$, if $\lambda+\rho$ is regular, there exists a unique $w\in W$ which makes $w(\lambda+\rho)$ dominant regular, and for this $w$, $w(\lambda+\rho)-\rho$ is a dominant weight.

\begin{thm} [Borel-Weil-Bott]
If $\lambda\in\Lambda$ and if $\lambda+\rho$ is singular, then
$$
H^p(X, {\cal O}({\cal L}_{\lambda})) =0\ \ \ \text{for all $p \in \BB Z$}.
$$
If $\lambda+\rho$ is regular, let $w$ be the unique element of $W$ such that $w(\lambda+\rho)$ is dominant, and define $p(\lambda) = \# \{ \alpha \in \Phi^+ \,\mid\, (\lambda+\rho, \alpha)<0 \}$. In this situation,
$$
H^p(X, {\cal O}({\cal L}_{\lambda})) =
\begin{cases}
\ \ \text{is non-zero, irreducible, of highest}
\\
\qquad\qquad\qquad\ \ \text{weight $w(\lambda +\rho)-\rho$}\ \ \  & \text{if $p = p(\lambda)$;} 
\\
\\
\ \ 0 & \text{if $p \ne p(\lambda)$.}
\end{cases}
$$
\end{thm}

The description of the highest weight as $w(\lambda +\rho)-\rho$ -- rather than $w\lambda$, for instance -- makes the statement compatible with Serre duality, as it has to be.

Bott proved this result by reducing it to the Borel-Weil theorem. The mechanism is a spectral sequence which relates the cohomology groups of line bundles ${\cal L}_{\lambda}$, ${\cal L}_{s(\lambda+\rho)-\rho}$, corresponding to parameters related by a so-called simple Weyl reflection $s\in W$ \cite{Bott}.  An outline of Bott's argument can be found in \cite{BSch}.

The Borel-Weil theorem alone suffices to realize all irreducible representations of $U_{\BB R}$. Bott's contribution is important for other reasons. At the time, it made it possible to compute some previously unknown cohomology groups of interest to
algebraic geometers. The Borel-Weil-Bott theorem made the flag varieties test cases for the fixed point formula for the index of elliptic operators. In representation theory, the theorem gave the first indication that higher cohomology groups might be useful in constructing representations geometrically. That turned out to be the case, for the representations of the discrete series of a non-compact reductive group, for example.

\begin{section}
{Representations of Reductive Lie Groups}
\end{section}

We now return to the case of a not-necessarily-compact, connected, linear reductive group $G_{\BB R}$. The notation of sections 1.2-3 applies. In particular, $K_{\BB R} \subset G_{\BB R}$ denotes a maximal compact subgroup, and $\g g$ is the complexified Lie algebra of $G_{\BB R}$.

\separate

\begin{subsection}
{Notions of Continuity and Admissibility, $K_{\BB R}$-finite and
$C^{\infty}$ Vectors}
\end{subsection}

Interesting representations of noncompact groups are typically infinite dimensional. To apply analytic and geometric methods, it is necessary to have a topology on the representation space and to impose an appropriate continuity condition on the representation in question. In the finite dimensional case, there is only one ``reasonable" topology and continuity hypothesis, but in the infinite dimensional case, choices must be made. One may want to study both complex and real representations. There is really no loss in treating only the complex case, since one can complexify  a real representation and regard the original space as an $\BB R$-linear subspace of its complexification. 

We shall consider representations on {\em complete locally convex
Hausdorff topological vector spaces over} $\BB C$. That includes complex Hilbert spaces, of course. Unitary representations are of particular interest, and one might think that Hilbert spaces constitute a large enough universe of representation spaces. It turns out, however, that even the study of unitary representations naturally leads to the consideration of other types of topological spaces, such as Fr\'echet spaces and DNF spaces. Most analytic arguments depend on completeness and the Hausdorff property. Local convexity is required to define the integral of vector-valued functions, which is a crucial tool in the study of representations of reductive groups -- see (\ref{operatorintegral}), for example.

Let $\operatorname{Aut}(V)$ denote the group of continuous, continuously invertible, linear maps from a complete locally convex Hausdorff topological vector space $V$ to itself; we do not yet specify a topology on this group. There are at least four reasonable notions of continuity one could impose on a homomorphism $G_{\BB R}\to \operatorname{Aut}(V)\,$:
\smallskip
\begin{itemize}
\item[\rm{a)}]
{\em continuity}: the action map $G_{\BB R} \times V \to V$ is continuous, relative to the product topology on $G_{\BB R} \times V$;
\item[\rm{b)}]
{\em strong continuity}: for every $v \in V$, $\,g \mapsto \pi(g)v\,$ is continuous as map from $G_{\BB R}$ to $V$;
\item[\rm{c)}]
{\em weak continuity}: for every $v \in V$ and every $\phi$ in the
continuous linear dual space $V'$, the complex-valued function $g \mapsto \langle \phi, \pi(g)v \rangle$ is continuous;
\item[\rm{d)}]
{\em continuity in the operator norm}, which makes sense only if $V$ is a Banach space; in that case, $\operatorname{Aut}(V)$ can be equipped with the norm topology, and continuity in the operator norm means that $\pi: G_{\BB R} \to \operatorname{Aut}(V)$ is a continuous homomorphism of topological groups.
\end{itemize}

\begin{rem}
The following implications hold for essentially formal reason:
\begin{equation*}
\text{ continuity } \Longrightarrow
\text{ strong continuity } \Longrightarrow
\text{ weak continuity}\,,
\end{equation*}
and if $V$ is a Banach space,
\begin{equation*}
\text{ continuity in the operator norm } \Longrightarrow
\text{ continuity}\,.
\end{equation*}
Also, if $V$ is a Banach space,
\begin{equation*}
\text{ continuity } \Longleftrightarrow
\text{ strong continuity } \Longleftrightarrow
\text{ weak continuity}\,.
\end{equation*}
In this chain of implications, {\it strong continuity}~$\Longrightarrow$~{\it continuity} follows relatively easily from the uniform boundedness principle, but the implication {\it weak continuity}~$\Longrightarrow$~{\it strong continuity}
is more subtle -- details can be found in {\rm \cite{War}}.
\end{rem}

\begin{ex}
The translation action of $(\BB R, +)$ on $L^p(\BB R)$ is continuous for $1\leq p < \infty$, but {\em not} continuous in the operator norm; for $p =\infty$ the translation action fails to be continuous, strongly continuous, even weakly continuous.
\end{ex}

Continuity in the operator norm is too much to ask for -- most of the representations of interest involve translation. Thus, from now on, ``representation" shall mean a continuous -- continuous in the sense described above -- linear action $\pi: G_{\BB R} \to \operatorname{Aut}(V)$ on a complete, locally convex Hausdorff space $V$. If $\pi$ is continuous, the dual linear action of the topological dual space $V'$, equipped with the strong dual topology\footnote{In the case of a Banach space, this is the dual Banach topology; for the general case, see \cite{Tr}, for example.}, need not be continuous. However, when $V$ is a reflexive Banach space, $V$ and $V'$ play symmetric roles in the definition of weak continuity; in this case, the dual action is also continuous, so there exists a ``dual representation" $\pi'$ of $G_{\BB R}$ on the dual Banach space $V'$.

An infinite dimensional representation $(\pi, V)$ typically has numerous invariant subspaces $V_1\subset V$, but the induced linear action of $G_{\BB R}$ on $V/V_1$ is a purely algebraic object unless $V/V_1$ is Hausdorff, i.e., unless $V_1\subset V$ is a closed subspace. For this reasons, the existence of a non-closed invariant subspace should not be regarded as an obstacle to irreducibility: $(\pi, V)$ is {\em irreducible} if $V$ has no proper {\em closed} $G_{\BB R}$-invariant subspaces. A representation $(\pi, V)$ has {\em finite length} if every increasing chain of closed $G_{\BB R}$-invariant subspaces breaks off after finitely many steps. One calls a representation $(\pi, V)$ {\em admissible} if $\dim_{\BB R} \Hom_{K_\BB R}(U,V) < \infty$ for every finite-dimensional irreducible representation $(\tau, U)$ of $K_{\BB R}$. Informally speaking, admissibility means that the restriction of $(\pi, V)$ to $K_{\BB R}$ contains any irreducible $K_{\BB R}$-representation only finitely often.

\begin{thm} [Harish-Chandra \cite{HC2}]
Every irreducible unitary representation $(\pi, V)$ of $G_{\BB R}$  is admissible.
\end{thm}

Harish-Chandra proved this theorem for a larger class of reductive Lie groups, not assuming linearity. Godement \cite{Go} gave a simplified, transparent argument for linear groups. Atiyah's lecture notes \cite{At} include Godement's argument and many related results.

Heuristically, admissible representations of finite length constitute the smallest class that is invariant under ``standard constructions" (in a very wide sense!) and contains the irreducible unitary representations. One should regard inadmissible irreducible representations as exotic. Indeed, the first example of an inadmissible irreducible representation on a Banach space -- a representation of the group $G_{\BB R}=SL(2,\BB R)$ -- is relatively recent \cite{Soe}, and depends on a counterexample to the invariant subspace problems for Banach spaces. All irreducible representations which have come up naturally in geometry, differential equations, physics, and number theory are admissible.

\begin{df}  \label{smooth}
Let $(\pi, V)$ be a representation of $G_{\BB R}$. A vector $v \in V$ is
\begin{itemize}
\item[{\rm a)}]
{\em $K_{\BB R}$-finite} if $v$ lies in a finite-dimensional
$K_{\BB R}$-invariant subspace;
\item[{\rm b)}]
a {\em $C^{\infty}$ vector} if $g \mapsto \pi(g)v$ is a $C^{\infty}$ map from $G_{\BB R}$ to $V$;
\item[{\rm c)}]
in the case of a Banach space $V$ only, an {\em analytic vector}, if $g \mapsto \pi(g)v$ is a $C^\omega$ map ($C^\omega$ means real analytic);
\item[{\rm d)}]
{a weakly analytic vector} if, for every $\phi \in V'$, the complex valued function $g \mapsto \langle \phi, \pi(g)v \rangle$ is
real analytic.
\end{itemize}
\end{df}

All reasonable notions of a real analytic $V$-valued map agree when $V$ is a Banach space, but not for other locally convex topological vector spaces. That is the reason for defining the notion of an analytic vector only in the Banach case. Surprisingly perhaps, even weakly real analytic functions with values in a Banach space are real analytic in the usual sense, i.e., locally representable by absolutely convergent vector valued power series -- see \cite[appendix]{La} for an efficient argument. In the Banach case, then, the notions of an analytic vector and of a weakly analytic coincide; for other representations, the former is not defined, but the latter still makes sense.

As a matter of self-explanatory notation, we write $V_{K_{\BB R}-\text{finite}}$ for the space of $K_{\BB R}$-finite vectors in $V$.

\begin{thm} [Harish-Chandra \cite{HC2}]\label{densityof}
If $(\pi, V)$ is an admissible representation,
\begin{itemize}
\item[{\rm a)}]
$V_{K_{\BB R}-\text{finite}}$ is a dense subspace of $V$;
\item[{\rm b)}]
every $v \in V_{K_{\BB R}-\text{finite}}$ is both a $C^{\infty}$ vector and a weakly analytic vector.
\end{itemize}
\end{thm}

Let us sketch the proof. For every $f\in C_c(G_{\BB R})=$ space compactly supported continuous functions on $G_{\BB R}$, the operator valued integral
\begin{equation}\label{operatorintegral}
\pi(f)\,\in\, \operatorname{End}(V)\,,\qquad \pi(f)v\ = \ \int_{G_{\BB R}} f(g)\,\pi(g)\,v\,dg \qquad(\,v\in V\,),
\end{equation} 
is well defined and convergent. If $f\in C_c^\infty(G_{\BB R})$, for any $v\in V$, $\pi(f)v$ is a $C^\infty$ vector, as was first observed by G\. arding. When $f=f_n$ runs through an ``approximate identity" in $C_c^\infty(G_{\BB R})$, the sequence $\pi(f_n)v$ converges to $v$. This much establishes G\. arding's theorem -- the space of $C^\infty$ vectors is dense in $V$. In analogy to the operators $\pi(f)$, one can also define $\pi(h)$ for $h\in C(K_{\BB R})$. Because of the Stone-Weierstrass theorem, the space of left $K_{\BB R}$-finite functions is dense in $C(K_{\BB R})$. Letting $h=h_n$ run through an approximate identity in $C(K_{\BB R})$, consisting of left $K_{\BB R}$-finite functions $h_n$, one can approximate any $v\in V$ by the sequence of vectors $\pi(h_n)v$, which are all $K_{\BB R}$-finite. This proves the density of the $K_{\BB R}$-finite vectors. One can combine this latter argument with G\.arding's, to approximate any $v\in V$ by a sequence of $K_{\BB R}$-finite $C^\infty$ vectors; when $v$ transforms according to a particular irreducible finite dimensional representation $\tau$ of $K_{\BB R}$, one can also make the approximating sequence lie in the space of $\tau$-isotypic vectors. This space is finite dimensional because of the admissibility hypothesis. The density of the subspace of $\tau$-isotypic $C^\infty$ vectors in the space of all $\tau$-isotypic vectors therefore implies that all $\tau$-isotypic vectors are $C^\infty$ vectors, for every $\tau$, so all $K_{\BB R}$-finite vectors are $C^\infty$ vectors. The functions $g \mapsto \langle \phi, \pi(g)v \rangle$, for $v\in V_{K_{\BB R}-\text{finite}}$ and $\phi\in V'$\!, satisfy elliptic differential equations with $C^\omega$ coefficients, which implies they are real analytic. Thus all $K_{\BB R}$-finite vectors are weakly analytic, as asserted by the theorem.

The theorem applies in particular to $K_{\BB R}$, considered as maximal compact subgroup of itself. Finite dimensional subspaces are automatically closed, so the density of $K_{\BB R}$-finite vectors forces any infinite dimensional representation $(\pi,V)$ of $K_{\BB R}$ to have proper closed invariant subspaces. In other words,

\begin{cor}
Every irreducible representation of $K_{\BB R}$ is finite dimensional. 
\end{cor}

\begin{subsection}
{Harish-Chandra Modules}
\end{subsection}

We continue with the notation of the previous section, but $(\pi, V)$ will now specifically denote an admissible representation of $G_{\BB R}$, and later an admissible representation of finite length. We write $V_{K_{\BB R}-\text{finite}}$ for the space of $K_{\BB R}$-finite vectors, and $V^\infty$ for the space of $C^\infty$ vectors. The latter is a $G_{\BB R}$-invariant subspace of $V$\!, which contains the former as $K_{\BB R}$-invariant subspace; both are dense in $V$. The Lie algebra $\g g_{\BB R}$ acts on $V^\infty$ by differentiation, and we extend this action by $\BB C$-linearity to the complexified Lie algebra $\g g$ on $V^\infty$~-- equivalently, $V^\infty$ has a natural structure of module over $U(\g g)$, the universal enveloping algebra of $\g g\,$. One might think that the $U(\g g)$-module $V^\infty$ is the right notion of infinitesimal representation attached to the global representation $(\pi, V)$. It has a very serious drawback, however: except in the finite dimensional case, $V^\infty$ may be highly reducible as $U(\g g)$-module even if $(\pi, V)$ is irreducible.

The first hint of a solution to this problem appeared in Bargmann's description of the irreducible unitary representations of $SL(2,\BB R)$ \cite{Ba}. Later formalized and developed by Harish-Chandra, it starts with the observation that
\begin{equation}
V_{K_{\BB R}-\text{finite}}\, \ \text{is a}\,\ U(\g g)-\text{submodule of}\,\ V^\infty\,.
\end{equation}
Indeed, the action map $U(\g g)\otimes V_{K_{\BB R}-\text{finite}} \to V^\infty$ is $K_{\BB R}$-invariant when $K_{\BB R}$ acts on $U(\g g)$ via the adjoint action and on $V$ and its subspaces via $\pi$. The image of the action map is therefore exhausted by finite dimensional, $K_{\BB R}$-invariant subspaces; in other words, the image of this action map lies in $V_{K_{\BB R}-\text{finite}}$.

Finite dimensional representations of the compact Lie group $K_{\BB R}$ extend naturally to its complexification (\ref{extensionofrep}). By definition, $V_{K_{\BB R}-\text{finite}}$ is the union of finite dimensional $K_{\BB R}$-invariant subspaces, so the $K_{\BB R}$-action on $V_{K_{\BB R}-\text{finite}}$ extends naturally to the complexification $K$ of $K_{\BB R}$. Even though $V_{K_{\BB R}-\text{finite}}$ has no natural Hausdorff topology~-- it is not closed in $V$ unless $\dim V < \infty$~-- it makes sense to say that $K$ acts holomorphically on $V_{K_{\BB R}-\text{finite}}\,$: like $K_{\BB R}$, $K$ acts {\em locally finitely}, in the sense that every vector lies in a finite dimensional invariant subspace; the invariant finite dimensional subspaces do carry natural Hausdorff topologies, and $K$ does act holomorphically on them. The Lie algebra $\g k$ has two natural actions on $V_{K_{\BB R}-\text{finite}}$, by differentiation of the $K$-action, and via the inclusion $\g k \subset \g g$ and the $U(\g g)$-module structure. These two actions coincide, essentially by construction. To simplify the notation, we denote the actions on $V_{K_{\BB R}-\text{finite}}$ by juxtaposition. With this convention,
\begin{equation}\label{adjoint}
k\,(\xi v)\ = \ (\operatorname{Ad}k\,\xi)(kv)\ \ \ \text{for all}\,\ \ k\in K\,,\ \ \xi \in U(\g g)\,,\ \ v\in V_{K_{\BB R}-\text{finite}}\,,
\end{equation}
as can be deduced from the well known formula $\exp(\operatorname{Ad}k\,\xi) = k \exp(\xi) k^{-1}$, for $\xi\in\g g_{\BB R}$, $k\in K_{\BB R}$. 
\medskip

\begin{df}\label{HCmodule}
A {\em $(\g g, K)$-module} is a complex vector space $M$, equipped with the structure of $U(\g g)$-module and with a linear action of $K$ such that:
\begin{itemize}
\item[{\rm a)}]
The action of $K$ is locally finite, i.e., every $m \in M$ lies in a finite dimensional $K$-invariant subspace on which $K$ acts holomorphically;
\item[{\rm b)}]
when the $K$-action is differentiated, the resulting action of $\g k$ agrees with the action of $\g k$ on $M$ via $\g k \hookrightarrow \g g$ and the $U(\g g)$-module structure.
\item[{\rm b)}]
the identity (\ref{adjoint}) holds for all $\,k\in K\,,\ \xi \in U(\g g)\,,\ v\in V_{K_{\BB R}-\text{finite}}$.
\end{itemize}
A {\em Harish-Chandra module} is a $(\g g, K)$-module $M$ which is  finitely generated over $U(\g g)$ and admissible, in the sense that every irreducible $K$-representation occurs in $M$ with finite multiplicity. 
\end{df}

The discussion leading up to the definition shows that the space of $K_{\BB R}$-finite vectors $V_{K_{\BB R}-\text{finite}}$ of an admissible representation $(\pi, V)$ is an admissible $(\g g,K)$-module. Very importantly,
\begin{equation}\label{HCdefined}
\begin{aligned}
&\text{the correspondence}\ \ \ \tilde V \ \mapsto \ \tilde V_{K_{\BB R}-\text{finite}} \ \ \ \text{sets up a bijection between closed}
\\
&\text{$G_{\BB R}$-invariant subspaces $\tilde V \subset V$ and $(\g g,K)$-submodules $\tilde V_{K_{\BB R}-\text{finite}} \subset V_{K_{\BB R}-\text{finite}}$}\,,
\end{aligned}
\end{equation}
as follows from the weakly analytic nature of $K_{\BB R}$-finite vectors (\ref{densityof}). When $(\pi, V)$ is not only admissible but also of finite length, every ascending chain of $(\g g,K)$-submodules of $V_{K_{\BB R}-\text{finite}}$ breaks off eventually. Because of the N\"otherian property of $U(\g g)$, that implies the finite gene\-ration of $V_{K_{\BB R}-\text{finite}}$ over $U(\g g)$. In short, the space of $K_{\BB R}$-finite vectors $V_{K_{\BB R}-\text{finite}}$ of an admissible representation of finite length is a Harish-Chandra module. The statement (\ref{HCdefined}) also implies that $(\pi, V)$ is irreducible if and only if $V_{K_{\BB R}-\text{finite}}$ is irreducible as $(\g g,K)$-module. This property of $V_{K_{\BB R}-\text{finite}}$ makes it the appropriate notion of the infinitesimal representation corres\-ponding to $(\pi, V)$.

From now on, we write $\HC(V)$ for the space of $K_{\BB R}$-finite vectors of an admissible representation $(\pi, V)$ and call $\HC(V)$ {\em the Harish-Chandra module of} $\pi$. The next statement formalizes the properties of Harish-Chandra modules we have mentioned so far: 

\begin{thm}[Harish-Chandra \cite{HC2}]\label{HCf} The association $V \mapsto \HC(V) = V_{K_{\BB R}-\text{finite}}$ establishes a covariant, exact, faithful functor
\begin{equation*}
\Bigl\{\begin{matrix}
\text{category of admissible $G_{\BB R}$-representations}
\\
\text{of finite length and $G_{\BB R}$-equivariant maps}
\end{matrix}\Bigr\} \ \overset{\HC}{\longrightarrow}   \
\Bigl\{\begin{matrix}
\text{category of Harish-Chandra modules} \\
\text{and $(\g g, K)$-equivariant linear maps}
\end{matrix}\Bigr\}.
\end{equation*}
\end{thm}

\begin{df}
Two admissible representations of finite length $(\pi_i,V_i)$, $i=1,2$,
are {\em infinitesimally equivalent} if $HC(V_1) \simeq HC(V_2)$.
\end{df}

Loosely speaking, infinitesimal equivalence means that the two representations are the same except for the choice of topology. A concrete example may be helpful. The group 
\begin{equation}\label{su1,1;1}
G_{\BB R} = SU(1,1) = \left\{ \left.
\begin{pmatrix} a & b \\ \bar b & \bar a \end{pmatrix} \ \right|\  a,b \in \BB C\,,\,\ |a|^2 - |b|^2 =1 \right\}
\end{equation}
has $G=SL(2,\BB C)$ as complexification, and is conjugate in $G$ to $SL(2,\BB R)$. As maximal compact subgroup, we choose the diagonal subgroup, in which case its complexification also consists of diagonal matrices:
\begin{equation}\label{su1,1;2}
K_{\BB R} = \left\{ \left. k_\theta = \begin{pmatrix} e^{i\theta} & 0 \\ 0 &  e^{-i\theta} \end{pmatrix} \ \right| \ \theta \in \BB R\,\right\}\ \cong \ U(1)\,,\ \ \ K = \left\{ \left.
\begin{pmatrix} a & 0 \\ 0 & a^{-1} \end{pmatrix} \ \right|\
a \in \BB C^*\,\right\}\ \cong \ \BB C^*\,.
\end{equation}
By linear fractional transformations, $SU(1,1)$ acts transitively on $D=$\ open unit disc in $\BB C\,$, with isotropy subgroup $K_{\BB R}$ at the origin. Left translation on $D\cong SU(1,1)/K_{\BB R}$ induces a linear action $\ell$ of $SU(1,1)$ on $C^\infty(D)$,
\begin{equation}\label{su1,1;3}
(\ell(g) f)(x) \ =_{\text{def}}\ f(g^{-1} \cdot z)\qquad (\,g \in SU(1,1)\,, \,\ f \in C^\infty(D)\,, \,\ z \in D\,),
\end{equation}
and on the subspace
\begin{equation}\label{su1,1;4}
H^2(D) \ =_{\text{def}}\  \text{space of holomorphic functions on $D$ with $L^2$ boundary values},
\end{equation}
topologized by the inclusion $H^2(D)\hookrightarrow L^2(S^1)$. One can show that both actions are representations, i.e., they are continuous with respect to the natural topologies of the two spaces.

Recall the definition of $k_{\theta}$ in (\ref{su1,1;2}). Since  $\ell(k_{\theta})z^n = e^{-2in\theta} z^n$, $\,f\in H^2(D)$ is $K_{\BB R}$-finite if and only if $f$ has a finite Taylor series at the origin, i.e., if and only if $f$ is a polynomial:
\begin{equation}\label{su1,1;5}
H^2(D)_{K_{\BB R}-\text{finite}} \ =\ \BB C [z]\,.
\end{equation}
In particular, $(\ell,H^2(D))$ is admissible. This representation is not irreducible, since $H^2(D)$ contains the constant functions $\BB C$ as an obviously closed invariant subspace. It does have finite length; in fact, the quotient $H^2(D)/\BB C$ is irreducible, as follows from a simple infinitesimal calculation in the Harish-Chandra module $\HC(H^2(D))=\BB C[z]$.

Besides $V=H^2(D)$, the action (\ref{su1,1;3}) on each of the
following spaces, equipped with the natural topology in each case, defines a representation of $SU(1,1)$:
\begin{itemize}
\item[{\rm a)}]
$H^p(D)$\ = space of holomorphic functions on $D$ with $L^p$ boundary values, $1 \le p \leq \infty$;
\item[{\rm b)}]
$H^\infty(D)$\ = space of holomorphic functions on $D$ with
$C^{\infty}$ boundary values;
\item[{\rm c)}]
$H^{-\infty}(D)$\ = space of holomorphic functions on $D$ with distribution boundary values;
\item[{\rm d)}]
$H^\omega(D)$\ = space of holomorphic functions on $D$ with real analytic
boundary values;
\item[{\rm e)}]
$H^{-\omega}(D)$\ = space of all holomorphic functions on $D$.
\end{itemize}
Taking boundary values, one obtains inclusions $H^p(D)\hookrightarrow L^p(S^1)$, which are equivariant with respect to the action of $SU(1,1)$ on $L^p(S^1)$ by linear fractional transformations. The latter fails to be continuous when $p=\infty$, but that is not the case for the image of $H^\infty(D)$ in $L^\infty(S^1)$. Essentially by definition, every holomorphic function on $D$ has hyperfunction boundary values. This justifies the notation $H^{-\omega}(D)$; the superscript $-\omega$ stands for hyperfunctions. One can show that $H^\infty(D)$ is the space of $C^\infty$ vectors for the Hilbert space representation $(\ell,H^2(D))$.

Arguing as in the case of $H^2(D)$, one finds that the representation $\ell$ of $SU(1,1)$ on each of the spaces a)-e) has $\BB C[z]$ as Harish-Chandra module, so all of them are infinitesimally equivalent. This is the typical situation, not just for $SU(1,1)$, but for all groups $G_{\BB R}$ of the type we are considering: every infinite dimensional admissible representation $(\pi,V)$ of finite length lies in an infinite family of representations, all infinitesimally equivalent, but pairwise non-isomorphic. In the context of unitary representations the situation is different:

\begin{thm} [Harish-Chandra \cite{HC2}]
If two irreducible unitary representations are infinitesimally equivalent, they are isomorphic as unitary representations.
\end{thm}

If we were dealing with finite dimensional representations, this would follow from an application of Schur's Lemma. Schur's lemma, it should be recalled, is a consequence of the existence of eigenvalues of endomorphisms of finite dimensional vector spaces over $\BB C$. In the setting of Harish-Chandra modules, endomorphisms are in particular $K_{\BB R}$-invariant, and must therefore preserve subspaces on which $K_{\BB R}$ acts according to any particular irreducible representation of $K_{\BB R}$. These subspaces are finite dimensional, and their direct sum, over all irreducible representations of $K_{\BB R}$ up to equivalence, is the Harish-Chandra module in question. Thus endomorphisms of Harish-Chandra modules can be put into Jordan canonical form, even though the modules are infinite dimensional. In short, there exists a version of Schur's lemma for Harish-Chandra modules; it plays the crucial role in the proof of the theorem.

According to results of Casselman \cite{Ca1}, every Harish-Chandra module $M$ has a {\em globalization}, meaning an admissible $G_{\BB R}$-representations $(\pi,V)$ of finite length, such that $\HC(V)=M$. That makes the following two problems equivalent:
\begin{itemize}
\item[{\rm a)}]
Classify irreducible admissible representations of $G_{\BB R}$, up to infinitesimal equivalence;
\item[{\rm b)}]
Classify irreducible Harish-Chandra modules for the pair $(\g g,K)$.
\end{itemize}
Under slight additional hypotheses, problem b) has been have been solved by Langlands \cite{Ll}, Vogan-Zuckerman \cite{Vo}, and Beilinson-Bernstein \cite{BB1,HMSW1}, by respectively analytic, algebraic, and geometric means. The three solutions give formally different answers, which can be related most easily in terms of the Beilinson-Bernstein construction \cite{HMSW2}. We shall discuss these matters in section 4.

A Harish-Chandra module corresponds to an irreducible unitary representation if and only if it carries an invariant, positive definite hermitian form -- invariant in the sense that the action of $K_{\BB R}$ preserves it, and that every $\xi\in\g g$ acts as a skew-hermitian transformation. Whether a given irreducible Harish-Chandra module carries an invariant, nontrivial, possibly indefinite hermitian form is easy to decide: the Harish-Chandra module needs to be conjugate-linearly isomorphic to its own dual. When a nontrivial invariant hermitian form exists, it is unique up to scaling. The classification of the irreducible unitary representations of $G_{\BB R}$ comes down to determining which invariant hermitian forms have a definite sign. Many results in this direction exist, most of them due to Vogan and Barbasch, but a general answer is not in sight, not even a good general conjecture.

We mentioned already that every Harish-Chandra module $M$ has a {\em globalization}. It is natural to ask if a globalization can be chosen in a functorial manner -- in other words, whether the functor $\HC$ in Theorem \ref{HCf} has a right inverse. Such functorial globalizations do exist. Four of them are of particular
interest, the $C^{\infty}$ and $C^{-\infty}$ globalizations
of Casselman-Wallach \cite{Ca2,Wal}, as well as the minimal globalization and the maximal globalization \cite{Sch,KSch}. All four are {\em topologically exact}, i.e., they map exact
sequences of Harish-Chandra modules into exact sequences of representations in which every morphism has {\em closed range}. The main technical obstacle in constructing the canonical globalizations is to establish this closed range property. In the case of an admissible representation $(\pi,V)$ of finite length, on a reflexive Banach space $V$, the $C^{\infty}$ globalization of $\HC(V)$ is topologically isomorphic to the space of $C^\infty$ vectors $V^\infty$. Similarly the minimal globalization is topologically isomorphic to the space of analytic vectors $V^\omega$; both have very naturally defined topologies. The other two constructions are dual to these: the $C^{-\infty}$ globalization is isomorphic to $((V')^\infty)'$, the strong dual of the space of $C^\infty$ vectors of the dual representation $(\pi',V')$, and the minimal globalization is similarly isomorphic to $((V')^\omega)'$. In the case of the earlier example of $(\ell,H^2(D))$, the four globalizations of $\HC(H^2(D))$ can be identified with $H^\infty(D)$, $H^\omega(D)$, $H^{-\infty}(D)$, and $H^{-\omega}(D)$, respectively.

\begin{section}
{Geometric Constructions of Representations}
\end{section}

\setcounter{equation}0
In this section we shall freely use the notational conventions of the preceding sections. To simplify the discussion, we suppose
\begin{equation}\label{connected}
\text{the complexification $G$ of $G_{\BB R}$ is connected}.
\end{equation}
That is the case for $G_{\BB R}=GL(n,\BB R)$, for example: the group itself is not connected, but it does have a connected complexification. The hypothesis (\ref{connected}) in particular implies:
\begin{equation}\label{connected'}
\text{for every $g\in G_{\BB R}$}\,,\ \ \Ad g \, : \, \g g \ \longrightarrow\ \g g\ \ \ \text{is an inner automorphism}.
\end{equation}
This latter condition is important; without it, irreducible Harish-Chandra modules need not have infinitesimal characters -- see definition (\ref{infcharacter}) below. From a technical point of view, the weaker condition (\ref{connected'}) suffices entirely for our purposes, at the cost of additional terminology and explanations. We assume (\ref{connected}) only to avoid these. The compact real form $U_{\BB R}\subset G$ is then also connected, as we had assumed in section 2.

Recall the notion of a Cartan subalgebra of the complex reductive Lie algebra $\g g\,$: a maximal abelian subalgebra $\g h\subset \g g$ such that $\Ad \xi : \g g \to \g g$ is diagonalizable, for every $\xi \in \g h$. Any two Cartan subalgebras are conjugate under the adjoint action of $G$. The complexified Lie algebra $\g t$ of a maximal torus $T_{\BB R} \subset U_{\BB R}$ is a particular example of a Cartan subalgebra. Since $G$ acts on the set of compact real forms by conjugation, every Cartan subalgebra of $\g g$ arises as the complexified Lie algebra of a maximal torus in some compact real form of $G$. In particular the discussion in sections 2.1-2 applies to any Cartan subalgebra. With $\g t$ and $T_{\BB R} \subset U_{\BB R}$ as above, there are two potential notions of Weyl group, namely the ``compact Weyl group" $W(U_{\BB R},T_{\BB R})=N_{U_{\BB R}}(T_{\BB R})/T_{\BB R}$ and the ``complex Weyl group" $W(G,T)=N_G(T)/T$, with $T=$ complexification of $T_{\BB R}$. They coincide, in fact, and we shall denote both by the symbol $W$.
\bigskip

\begin{subsection}
{The Universal Cartan Algebra and Infinitesimal Characters}
\end{subsection}

By definition, the flag variety $X$ of the complex reductive Lie algebra $\g g$ parameterizes the Borel subalgebras of $\g g\,$:
\begin{equation}
X \ni x \quad \longleftrightarrow \quad \g b_x \subset \g g\,.
\end{equation}
Define
\begin{equation}
\g h_x = \g b_x / [\g b_x, \g b_x].
\end{equation}
This quotient is independent of $x$ in the following equivalent
senses:
\medskip
\begin{itemize}
\item[{\rm a)}]
if $g\cdot x = y$ for some $g \in G$ and $x, y \in X$,
the map $\g h_x \to \g h_y$ induced by $\Ad g :\g b_x \to\g b_y$ depends only on $x$ and $y$, not on the particular choice of $g\,$;
\item[{\rm b)}]
$\g h_x$ is the fiber at $x$ of a canonically flat
holomorphic vector bundle over $X$;
\item[{\rm c)}]
let $\g t \subset \g g$ be a Cartan subalgebra and $\Phi^+ \subset \Phi(\g g,\g t)$ a positive root system. Then $\g b_0 = \g t \oplus
\left( \oplus_{\alpha \in \Phi^+}\, \g g^{-\alpha} \right)$ is a Borel subalgebra, with $[\g b_0,\g b_0]=\oplus_{\alpha \in \Phi^+} \g\, g^{-\alpha}$. The resulting isomor\-phism $\g h_0=\g b_0/[\g b_0,\g b_0]\cong \g t$ depends only on the choice of $\Phi^+$.
\end{itemize}

We write $\g h$ instead of $\g h_x$ to signify independence of $x$.
This is the {\em universal Cartan algebra}. It is {\em not a
subalgebra} of $\g g$, but $\g h$ {\em is} canonically isomorphic to any {\em ordered Cartan subalgebra}, i.e., to any Cartan subalgebra $\g t \subset \g g$ with a specified choice of positive root system $\Phi^+$. We use the canonical isomorphism between $\g h$ and any ordered Cartan subalgebra $(\g t,\Phi^+)$ to transfer from $\g t$ to $\g h$ the weight lattice, root system, positive root system, and Weyl group. In this way we get the {\em universal weight lattice} $\Lambda \subset \g h^*$, the {\em universal root system} $\Phi \subset \Lambda$, the {\em universal positive root system} $\Phi^+ \subset \Phi$, and the {\em universal Weyl group} $W$, which acts on $\g h$ and dually on $\g h^*$, leaving invariant both $\Lambda$ and $\Phi$. Moreover, there exists a $W$-invariant, positive inner product $(.,.)$ on the $\BB R$-linear subspace $\BB R\otimes_{\BB Z}\Lambda$ which depends on the choice of $S$ in (\ref{Sdef}), but on nothing else. Going back to (\ref{Llambda}), we see that the parametrization of the $G$-equivariant holomorphic line bundles $\mathcal L_\lambda \to X$ in terms of $\lambda$ becomes completely canonical when we regard $\lambda$ as lying in universal weight lattice $\Lambda$. 

We can use these ideas to characterize the so-called Harish-Chandra isomorphism. Let $Z(\g g)$ denote the center of the universal enveloping algebra $U(\g g)$, and $S(\g h)^W$ the algebra of $W$-invariants in the symmetric algebra of $\g h$, or what comes to the same, the algebra of $W$-invariant polynomial functions on $\g h^*$. By differentiation of the $G$-action, $\g g$ acts on holomorphic sections of $\mathcal L_\lambda$ as a Lie algebra of vector fields. This induces an action of $U(\g g)$, and therefore also of $Z(\g g)$, on the sheaf of holomorphic sections ${\mathcal O} ({\mathcal L}_{\lambda})$.

\begin{thm} [Harish-Chandra \cite{HC1}]\label{defofgamma}
There exists a canonical isomorphism
$$
\gamma: Z(\g g) \ \overset{\sim}{\longrightarrow} \ S(\g h)^W
$$
such that, for any $\lambda \in \Lambda$, any $\zeta \in Z(\g g)$ acts on the sheaf of holomorphic sections ${\cal O} ({\cal L}_{\lambda})$ as multiplication by the scalar $\gamma(\zeta)(\lambda+\rho)$.
\end{thm}

In this statement, $\gamma(\zeta)(\lambda+\rho)$ refers to the value of the $W$-invariant polynomial function $\gamma(\zeta)$ at the point $\lambda+\rho$. More generally this makes sense for elements of $\g h^*$\,: every $\lambda \in \g h^*$ determines a character\footnote{in the present context, ``character" means algebra homomorphism into the one dimensional algebra $\BB C$.}
\begin{equation}
\chi_{\lambda}\, :\, Z(\g g) \ \longrightarrow\  \BB C\,,
\qquad \chi_{\lambda}(\zeta) = \gamma(z)(\zeta).
\end{equation}
In view of Harish-Chandra's theorem, $\chi_{\lambda} = \chi_{\mu}$ if and only if $\lambda = w \cdot \mu$ for some $w \in W$, and every character of $Z(\g g)$ is of this type.

\begin{df}\label{infcharacter}
One says that a Harish-Chandra module $M$ has an {\em infinitesimal character} if $Z(\g g)$ acts on $M$ via a character.
\end{df}

When the infinitesimal character exists, it can of course be expressed as $\chi_\lambda$, for some $\lambda\in \Lambda$. Applying Schur's lemma for Harish-Chandra modules, as was explained in section 3.2, one finds that every irreducible Harish-Chandra
module does have an infinitesimal character.

\begin{subsection}
{Twisted ${\cal D}$-modules}
\end{subsection}

The flag variety $X$ is projective, and thus in particular has an algebraic structure. The complex linear reductive group $G$ has an algebraic structure as well, and the action of $G$ on $X$ is algebraic. The $G$-equivariant line bundles $\mathcal L_\lambda \to X$ are associated to algebraic characters of the structure group of the algebraic principal bundle $G \to G/B$, so they, too, have algebraic structures, and $G$ acts algebraically also on these line bundles. In the following, we equip $X$ with the Zariski topology,
and all sheaves are understood to be sheaves relative to the Zariski topology. In the current setting, ${\cal O}$ denotes the sheaf of algebraic functions on $X$ and $\mathcal O(\mathcal L_\lambda)$ the sheaf of algebraic sections of $\mathcal L_\lambda$.

The locally defined linear differential operators on $X$ with
algebraic coefficients constitute a sheaf of algebras, customarily denoted by  ${\cal D}$. By definition, $\mathcal D$ acts on ${\cal O}$; in more formal language, $\mathcal O$ is a sheaf of modules over the sheaf of algebras $\mathcal D$. For $\lambda \in \Lambda$,
\begin{equation}\label{twisting}
{\cal D}_{\lambda}\ =\ {\cal O}({\cal L}_{\lambda}) \otimes_{\cal O}
{\cal D} \otimes_{\cal O} {\cal O}({\cal L}_{-\lambda})
\end{equation}
is also a sheaf of algebras, the sheaf of linear differential operators acting on sections of ${\cal L}_{\lambda}$. It is a so-called {\em twisted sheaf of differential operators}. Since we can think of ${\cal O}$ as differential operators of degree zero,
there exists a natural inclusion ${\cal O} \hookrightarrow {\cal D}_{\lambda}$. By differentiation of the $G$-action, every $\xi\in\g g$ determines a globally defined first order differential operator acting on $\mathcal O({\cal L}_{\lambda})$. In this way we get a canonical homomorphism of Lie algebras
\begin{equation}\label{gtogammaD}
\g g\ \longrightarrow\ \Gamma {\cal D}_{\lambda}\ =\ H^0(X,{\cal D}_{\lambda})\,;
\end{equation}
here, as usual, we give the associative algebra $\Gamma {\cal D}_{\lambda}$ the additional structure of a Lie algebra by taking commutators of differential operators. When $\g g$ is semisimple, this morphism is injective. In any case, it induces
\begin{equation}\label{UgtogammaD}
U(\g g)\ \longrightarrow\ \Gamma {\cal D}_{\lambda}\ =\ H^0(X,{\cal D}_{\lambda})\,,
\end{equation}
a homomorphism of associative algebras.

Until now we have supposed that $\lambda$ lies in the weight lattice. However, there is a natural way to make sense of the sheaf of algebras ${\cal D}_{\lambda}$ any $\lambda \in \g h^*$. In terms of local coordinates, the twisting operation (\ref{twisting}) involves taking logarithmic derivatives, so the lattice parameter $\lambda$ occurs polynomially; one can therefore let it take values in $\g h^*$. With this extended definition, the natural inclusion ${\cal O} \hookrightarrow {\cal D}_{\lambda}$ and the morphism (\ref{UgtogammaD}) exist just as before. 
\bigskip

\begin{thm} [Beilinson-Bernstein \cite{BB1}]\label{gammaofdlambda}   For any $\lambda \in \g h^*\!$, the morphism (\ref{UgtogammaD}) induces
$$
U_{\lambda+\rho} \ =_{\text{def}}\ U(\g g) \,/ \,\text{ideal generated by $\{\zeta - \chi_{\lambda+\rho}(\zeta) \mid \zeta \in Z(\g g) \}$}\ \ \overset{\sim}{\longrightarrow}\ \ \Gamma{\cal D}_{\lambda}\,. 
$$
The higher cohomology groups of $\,{\cal D}_{\lambda}$ vanish:  $H^p(X,{\cal D}_{\lambda})=0$ for $p>0$\,.
\end{thm}

Note that any $U_{\lambda+\rho}$-module can be regarded as a $U(\g g)$-module with infinitesimal character $\chi_{\lambda+\rho}$, and vice versa. Following the usual custom, we shall use the terminology ``$\mathcal D_\lambda$-module" as shorthand for ``sheaf of $\mathcal D_\lambda$-modules". A $\mathcal D_\lambda$-module is said to be {\em coherent} if it is coherent over the sheaf of algebras $\mathcal D_\lambda$ -- in other words, if locally around any point, it can be presented as the quotient of a free $\mathcal D_\lambda$-module of finite rank, modulo the image of some other free $\mathcal D_\lambda$-module of finite rank. Theorem \ref{gammaofdlambda} makes it possible to define the two functors
\begin{equation}\label{equivofcat}
\left\{ \begin{matrix}\text{category of finitely} \\ \text{ generated $U_{\lambda+\rho}$-modules}\end{matrix} \right\}\ \ \overset{\displaystyle{\overset{\Delta}{\longrightarrow}}} {\underset{\Gamma}{\longleftarrow}}\ \ 
\left\{ \begin{matrix}
\text{category of} \\
\text{coherent ${\cal D}_{\lambda}$-modules} \end{matrix} \right\}\,,
\end{equation}
with $\Gamma = H^0(X,\,\cdot\,)=$ global sections functor, and
\begin{equation}\label{localization}
\Delta M\ \ = \ \ {\cal D}_{\lambda} \otimes_{U_{\lambda+\rho}} M\,;
\end{equation}
$\Delta$ is called the ``localization functor". Because of (\ref{defofgamma}), the category on the left in (\ref{equivofcat}) depends only on the $W$-orbit of $\lambda+\rho\,$; the category on the right, on the other hand, depends on $\lambda$ itself.

Recall the definition of the $W$-invariant inner product $(\cdot,\cdot)$ on the real form $\BB R \otimes_{\BB Z}\Lambda\subset\g h^*$ in section 4.1. We shall use the same notation to denote the bilinear -- not hermitian! -- extension of the inner product to the complex vector space $\g h^*$. In analogy to our earlier terminology we call $\lambda\in\g h^*$ {\em regular} if $(\lambda,\alpha)\neq 0$ for all $\alpha\in\Phi$, and otherwise {\em singular}.

\begin{df}\label{int_dom}
An element $\lambda \in \g h^*$ is {\em integrally dominant} if\ \ $
2 \frac {(\lambda, \alpha)}{(\alpha,\alpha)} \notin \BB Z_{<0}$\,\ for all $\,\alpha \in \Phi^+$.
\end{df}

\begin{rem}
For $\lambda \in \Lambda$ these quotients are integers, and for a generic $\lambda\in\g h^*$ all of them are non-integral. In every case there exists $\,w\in W$ such that $w(\lambda+\rho)$ is integrally dominant.
\end{rem}

\begin{thm} [Beilinson-Bernstein \cite{BB1}]\label{AB}
Let ${\cal S}$ be a coherent ${\cal D}_{\lambda}$-module.
\begin{itemize}
\item[{\rm A)}] If $\lambda+\rho$ is integrally dominant and regular, the global sections of ${\cal S}$ generate its stalks;
\item[{\rm B)}] If $\lambda+\rho$ is integrally dominant, then
$H^p(X,{\cal S}) =0$ for all $p\ne 0$.
\end{itemize}
\end{thm}

The conclusion of A) means that the stalk $\mathcal S_x$ at any $x\in X$ is generated over the ring $(\mathcal D_\lambda)_x$ by the image of $\Gamma \mathcal S$ in $\mathcal S_x$. Note the formal analogy with Cartan's theorems A and B for coherent analytic sheaves on Stein manifolds -- see \cite{GuRo}, for example.  The next statement follows quite directly from the theorem:

\begin{cor}\label{ABcor}
If $\lambda+\rho$ is integrally dominant and regular, the localization functor $\Delta$ defines an equivalence of categories, with inverse functor $\Gamma$.
\end{cor}

The hypothesis of regularity is essential. For example, when the canonical bundle of $X$ has a $G$-equivariant square root $\mathcal L_{-\rho}$, the Borel-Weil-Bott theorem asserts that the $\mathcal D_{-\rho}$-module ${\cal O}({\cal L}_{-\rho})$ has no non-zero global sections, nor even non-zero higher cohomology, but $-\rho+\rho = 0$ {\em is} integrally dominant. A refined version of the corollary does apply whenever $\lambda+\rho$ is integrally dominant but singular. The idea is to construct a quotient category of the category of coherent $\mathcal D_\lambda$-modules, by dividing out the subcategory of sheaves with trivial cohomology. The greater subtlety of this situation merely reflects, and even explains, a familiar fact: in the study of $U(\g g)$-modules with singular infinitesimal character, one encounters difficulties not present in the regular case. The category of coherent $\mathcal D_\lambda$-modules ``knows nothing" about regularity or singularity; for any $\mu\in\Lambda$,
\begin{equation}
\left\{ \begin{matrix}
\text{category of} \\
\text{coherent ${\cal D}_{\lambda}$-modules} \end{matrix} \right\}\ \ni \ \mathcal S\ \mapsto\ \mathcal O(\mathcal L_\mu)\otimes_{\mathcal O}\mathcal S\ \in \ \left\{ \begin{matrix}
\text{category of} \\
\text{coherent ${\cal D}_{\lambda+\mu}$-modules} \end{matrix} \right\}
\end{equation}
defines an equivalence of categories, and $\mu$ can always be chosen so as to make $\lambda+\mu+\rho$ integrally dominant and regular. The existence of sheaves without cohomology carries sole responsibility for the more complicated nature of the singular case!

The corollary can be extended in another direction. When $\lambda+\rho$ fails to be integrally dominant, the localization functor $\Delta$ becomes an equivalence of categories on the level of derived categories. Heuristically, and quite imprecisely, that means replacing $\Gamma{\cal S}$ by the formal Euler characteristic $\sum_p (-1)^p H^p (X,{\cal S})$. In other words, the derived category of coherent $\mathcal D_{\lambda}$-modules is equivalent to the derived category of coherent $\mathcal D_{w(\lambda+\rho)-\rho}$-modules, for any $w\in W$. This latter equivalence can be described geometrically, in terms of Beilinson-Bernstein's {\em intertwining functors} \cite{BB2}. 

The corollary and the first of the two extensions makes it possible to translate problems about finitely generated $U(\g g)$-modules with an infinitesimal character into problems in algebraic geometry. What are advantages to working on the geometric side, rather than directly with $U(\g g)$-modules? For those of us who think geometrically, the geometric arguments seem far more transparent than their algebraic counterparts. More importantly, some results on $U(\g g)$-modules, which seem inaccessible by algebra, have been proved geometrically -- the proofs of the Kazhdan-Lusztig conjectures \cite{BB1,BK} and of the Barbasch-Vogan conjecture \cite{SchV} are particular examples.

\begin{subsection}
{Construction of Harish-Chandra Modules}  \label{HCmods}
\end{subsection}

The equivalence of categories (\ref{ABcor}) persists when certain additional ingredients are fed in on both sides. In the case of a Harish-Chandra module $M$ with infinitesimal character $\chi_{\lambda+\rho}$, the algebraic action\footnote{All finite dimensional representations of a complex linear reductive group are algebraic.} of the group $K$ on $M$ induces an algebraic action of $K$ on the $\mathcal D_\lambda$-module $\Delta M$. One might think that the admissibility of $M$ puts an additional restriction on $\Delta M$, but that is not the case: any $(\g g,K)$-module with an infinitesimal character is automatically admissible. By a {\em $({\cal D}_{\lambda}, K)$-module} one means a sheaf of ${\cal D}_{\lambda}$-modules, equipped with an algebraic action of $K$ that is compatible with the ${\cal D}_{\lambda}$-structure. Compatibility in the current setting is entirely analogous to the earlier notion of compatibility in the context of $(\g g,K)$-module. More precisely, the Lie algebra $\g k$ acts on $M$ both via the differentiation of the $K$-action and via (\ref{gtogammaD}) and the inclusion $\g k \hookrightarrow\g g$. These two actions of $\g k$ must agree, and the analogue of (\ref{adjoint}) must also be satisfied. With these additional ingredients, the functors $\Delta$ and $\Gamma$ in (\ref{equivofcat}) induce
\begin{equation}\label{equivofcat'}
\left\{ \begin{matrix}\text{category of Harish-Chandra modules} \\ \text{with infinitesimal character $\chi_{\lambda+\rho}$}\end{matrix} \right\}\ \ \overset{\displaystyle{\overset{\Delta}{\longrightarrow}}} {\underset{\Gamma}{\longleftarrow}}\ \ 
\left\{ \begin{matrix}
\text{category of coherent} \\
\text{$({\cal D}_{\lambda},K)$-modules} \end{matrix} \right\}\,.
\end{equation}
Corollary \ref{ABcor} has a counterpart for these restricted functors:
\medskip

\begin{cor}\label{ABcor'}
If $\lambda+\rho$ is integrally dominant and regular, the functor $\Delta$ in (\ref{equivofcat'}) defines an equivalence of categories, with inverse functor $\Gamma$.
\end{cor}

This corollary, too, has refinements that apply when either or both of the hypotheses of integral dominance and regularity are dropped. Under an equivalence of categories, irreducible objects correspond to irreducible objects, hence

\begin{cor}  \label{ABcorcor}
If $\lambda+\rho$ is integrally dominant and regular, the functor $\Delta$ establishes a bijection
$$
\Bigl\{\begin{matrix} \text{irreducible Harish-Chandra modules} \\
\text{with infinitesimal character $\chi_{\lambda+\rho}$}\end{matrix}\Bigr\}\ \ 
\underset{\displaystyle{\sim}}{\overset{\Delta}{\longrightarrow}}\ \ \Bigl\{\begin{matrix}
\text{irreducible} \\
\text{$({\cal D}_{\lambda}, K)$-modules}
\end{matrix}\Bigr\}.
$$
If $\lambda+\rho$ is integrally dominant but singular, $\Delta$ sets up a bijection between the set of all irreducible Harish-Chandra modules with infinitesimal character $\chi_{\lambda+\rho}$ on the one hand, and the set of irreducible $({\cal D}_{\lambda},K)$-modules with non-zero cohomology on the other.
\end{cor}

This latter corollary describes the irreducible Harish-Chandra modules in terms of irreducible $({\cal D}_{\lambda}, K)$-modules. But how does one construct such sheaves? Two properties of the $K$-action on $X$ simplify the problem. First of all,
\begin{equation}
\text{$K$ acts on $X$ with finitely many orbits.}
\end{equation}
Since $K$ acts algebraically, these orbits are algebraic subvarieties. Moreover,
\begin{equation}\label{affineemb}
\text{all $K$-orbits in $X$ are affinely embedded,}
\end{equation}
which means that they intersect any open affine subset $U \subset X$ in an affine set. As a negative example, we mention $\BB {CP}^n - \{\text{point}\}$, with $n\geq 2$, which is not affinely embedded in $\BB {CP}^n$. Such sets do arise as $K$-orbits in generalized flag varieties.

We now consider a particular irreducible $({\cal D}_{\lambda}, K)$-module $\mathcal S$. For geometric reasons, the support of ${\cal S}$ must consist of the closure of a single $K$-orbit $Q$. We let $\partial Q$ denote the boundary of $Q$ in $X$, i.e., $\partial Q=$ (closure of $Q)\,-Q$. The inclusion $j: Q \hookrightarrow X$ factors as a product
\begin{equation}
j \ = \ j_o \circ j_c\,,\ \ \ \text{with}\ \ j_c\,:\, Q \ \hookrightarrow \ X - \partial Q\ \ \ \text{and}\ \ j_o\,:\, X - \partial Q\ \hookrightarrow\ X\,;
\end{equation}
note that $j_c$ is a smooth closed embedding and $j_0$ an open embedding. Since ${\cal S}|_{X - \partial Q}$ is $K$-equivariant, irreducible, and supported on the $K$-orbit $Q$, ${\cal L}_{\lambda}$ must exist at least as a $K$-equivariant line bundle on the formal neighborhood of $Q$ in $X - \partial Q$, even if it does not exist as line bundle on all of $X$, and ${\cal S}|_{X - \partial Q}$ must be the $\mathcal D_\lambda$-module direct image of ${\cal O}_Q ({\cal L}_{\lambda}|_Q)$ under $j_c$ -- in formal notation,
\begin{equation}
{\cal S}|_{X - \partial Q}\ = \ j_{c+} {\cal O}_Q ({\cal L}_{\lambda}|_Q)\,.
\end{equation}
Depending on the orbit $Q$, this forces certain integrality conditions on $\lambda$; if these integrality conditions do not hold, no irreducible $({\cal D}_{\lambda}, K)$-module can have the closure of $Q$ as support. The ${\cal D}_{\lambda}$-module direct image $j_{c+} {\cal O}_Q ({\cal L}_{\lambda}|_Q)$ is easy to describe because $j_c$ is a smooth closed embedding: its sections can be expressed as normal derivatives, of any order, applied to sections of ${\cal L}_{\lambda}|_Q$ over $Q$. Since $\mathcal S$ is irreducible,
\begin{equation}
{\cal S}\ \hookrightarrow \ j_{o+} ({\cal S}|_{X \setminus \partial Q})\ = \ j_{o+} \circ j_{c+} {\cal O}_Q ({\cal L}_{\lambda}|_Q)\  = \ j_+ {\cal O}_Q ({\cal L}_{\lambda}|_Q)\,.
\end{equation}
In general, constructing the $\mathcal D$-module direct image under an open embedding requires passage to the derived category. In our situation, because of (\ref{affineemb}), the direct image exists as a bona fide $\mathcal D_\lambda$-module. One calls $j_+ {\cal O}_Q ({\cal L}_{\lambda}|_Q)$ the {\em standard sheaf} corresponding to the orbit $Q$, the parameter $\lambda\in\g h^*$, and one other simple datum that is necessary to pin down the meaning of ${\cal L}_{\lambda}|_Q$. The steps we outlined exhibit the irreducible $({\cal D}_{\lambda}, K)$-module $\mathcal S$ as the {\em unique irreducible subsheaf} of the standard sheaf $j_+ {\cal O}_Q ({\cal L}_{\lambda}|_Q)$.

At one extreme, the irreducible subsheaf ${\cal S}$ may coincide with the standard sheaf $j_+ {\cal O}_Q ({\cal L}_{\lambda}|_Q)$
in which it lies, and at the opposite extreme, it may be much smaller. This phenomenon is governed by the behavior of sections near $\partial Q$. Very roughly, if $j_+ {\cal O}_Q ({\cal L}_{\lambda}|_Q)$ has sections with various degrees of regularity along $\partial Q$, the unique irreducible subsheaf $\mathcal S$ consists of the ``most regular" sections; when all sections have the same degrees of regularity along $\partial Q$, the standard sheaf $j_+ {\cal O}_Q ({\cal L}_{\lambda}|_Q)$ is irreducible, hence equal to its unique irreducible subsheaf ${\cal S}$. The standard sheaf is more tractable than $\mathcal S$. In the crucial situation, when $\lambda + \rho$ is integrally dominant, one understands $H^0(X,j_+ {\cal O}_Q ({\cal L}_{\lambda}|_Q))$, the {\em standard module} corresponding to the given set of data, quite well \cite{HMSW1}. The space of sections $H^0(X,\mathcal S)$ is the unique irreducible submodule of the standard module -- or zero, which can happen only when $\lambda + \rho$ is singular. In view of (\ref{ABcorcor}), these results constitute a classification of the irreducible Harish-Chandra modules with infinitesimal character $\chi_{\lambda+\rho}$ -- the {\em Beilinson-Bernstein classification}.

Two other classification schemes, due to Landglands \cite{Ll} and Vogan-Zuckerman \cite{Vo}, predate Beilinson-Bernstein's. They, too, exhibit the irreducible Harish-Chandra modules as unique irreducible submodules, or dually as unique irreducible quotients, of certain standard modules. All three classifications obviously describe the same class of objects, but it is not clear a priori that the three types of standard modules agree or are dual to each other. That can be shown most transparently by geometric arguments \cite{HMSW2}. Depending on the interplay between the orbit $Q$ and the parameter $\lambda$, the higher cohomology groups of the standard sheaf $j_+ {\cal O}_Q ({\cal L}_{\lambda}|_Q)$ may vanish even when part B) of Theorem (\ref{AB}) does not apply directly. One extreme case, with ``$Q$ as affine as possible'', leads to the Langlands classification, the other, with ``$Q$ as close to projective as possible'', to Vogan-Zuckerman's. The Beilinson-Bernstein situation lies between these two, and all three can be related via the intertwining functors we mentioned earlier.

\begin{subsection}
{Construction of $G_{\BB R}$-representations}  \label{G_Rreps}
\end{subsection}

Just as one can attach Harish-Chandra modules to $K$-orbits in the flag variety, $G_{\BB R}$-represen\-tations arise from $G_{\BB R}$-orbits. There are finitely many such orbits, and they are real algebraic subvarieties. We now equip $X$ with the usual Hausdorff topology -- not the Zariski topology, as in the previous section.

To motivate the discussion, we first look at two special cases. At one extreme, let us consider a group $G_{\BB R}$, subject to the usual hypotheses and the condition (\ref{connected}). We suppose that $G_{\BB R}$ contains a compact Cartan subgroup, and we fix an open $G_{\BB R}$-orbit $S\subset X$. As subgroup of $G$, $G_{\BB R}$ acts on $\mathcal O(\mathcal L_\lambda)$, the sheaf of holomorphic sections of a $G$-equivariant line bundle $\mathcal L_\lambda$. This action induces a natural linear action on the cohomology groups $H^p(S,\mathcal O(\mathcal L_\lambda))$ over the open $G_{\BB R}$-orbit $S$. The cohomology can be computed from the complex of $\mathcal L_\lambda$-valued Dolbeault forms. It is far from obvious, but the coboundary operator $\bar\partial$ has closed range, giving the cohomology groups natural Fr\'echet topologies, with respect to which $G_{\BB R}$ acts continuously. The resulting representations are admissible, of finite length, with infinitesimal character $\chi_{\lambda+\rho}$. When $\lambda+\rho$ is antidominant regular\footnote{i.e., $(\lambda+\rho,\alpha)<0$ for all $\alpha\in\Phi^+$.}, the cohomology vanishes except in a single degree $s>0$, and in that degree $p=s$, $H^s(S,\mathcal O(\mathcal L_\lambda))$ is a discrete series representation -- more precisely, it is the maximal globalization of the Harish-Chandra module $\HC(V_{\lambda+\rho})$ of a discrete series representation $V_{\lambda+\rho}$. This construction provides a geometric realization, analogous to the Borel-Weyl-Bott theorem, of the entire discrete series. The reader may consult \cite{BSch} for references and further details.

At the other extreme, we suppose that $S\subset X$ is a closed orbit of a split group $G_{\BB R}$; here ``split" means that $G_{\BB R}$ contains a Cartan subgroup $A_{\BB R}$ such that every $a\in A_{\BB R}^0$ acts with real eigenvalues. Parenthetically we should remark that noncompact reductive groups $G_{\BB R}$ have finitely many Cartan subgroups, and may contain both a compact and a split Cartan subgroup. In the split case, the closed $G_{\BB R}$-orbit $S$ is necessarily a real form of the flag variety, i.e., a submanifold such that the holomorphic tangent space $T_xX$ at any $x\in S$ contains the tangent space $T_xS$ of $S$ as a real form: $T_xX = \BB C \otimes_{\BB R} T_xS$. For any $\lambda\in \Lambda$, $G_{\BB R}$ acts on $C^{-\omega}(S,\mathcal L_\lambda)$, the space of hyperfunction sections of $\mathcal L_\lambda$ over the real analytic manifold $S$. The hyperfunctions on a compact real analytic manifold carry a natural Fr\'echet topology. With respect to this topology, $G_{\BB R}$ acts continuously on $C^{-\omega}(S,\mathcal L_\lambda)$. The resulting representation is admissible, of finite length, with infinitesimal character $\chi_{\lambda+\rho}$; it belongs to the -- non-unitary, in general -- principal series of $G_{\BB R}$. Since $\lambda$ has been confined to the lattice $\Lambda$, only some principal series representations can be obtained this way. One gets the others by letting $\lambda$ range over $\g h^*$, the dual space of the universal Cartan, in which case $\mathcal L_\lambda$ still exists as $G_{\BB R}$-equivariant real analytic line bundle over the closed orbit $S$. However, for the moment we still want to suppose $\lambda\in\Lambda$, so that $\mathcal L_\lambda$ is well defined even as $G$-equivariant holomorphic line bundle on $X$.

At first glance, the cohomology groups $H^p(S,\mathcal O(\mathcal L_\lambda))$ over an open $G_{\BB R}$-orbit $S$ and the space of hyperfunction sections $C^{-\omega}(S,\mathcal L_\lambda)$ over a closed orbit of a split group $G_{\BB R}$ might not seem to fit easily into a common framework. However, both can be expressed as $\operatorname{Ext}$ groups,
\begin{equation}\label{extconstruction}
\begin{aligned}
H^p(S,\mathcal O(\mathcal L_\lambda))\ &= \, \operatorname{Ext}^p(j_!\BB C_S,\mathcal O(\mathcal L_\lambda))\ \ \text{if $S$ is an open $G_{\BB R}$-orbit},
\\
C^{-\omega}(S,\mathcal L_\lambda)\ &= \ \operatorname{Ext}^n(j_!\BB C_S,\mathcal O(\mathcal L_\lambda))\ \ \text{if $S$ is closed,  $G_{\BB R}$ split, and $n=\dim_{\BB C}X$};
\end{aligned}
\end{equation} 
in both cases $j$ denotes the inclusion $S\hookrightarrow X$, $\,\BB C_S$ the constant sheaf on $S$ with fiber $\BB C$, and $j_!\BB C_S$ the sheaf on $X$ obtained by taking the direct image with proper supports. One can describe the $\operatorname{Ext}$ groups equivalently as the right derived functors of $\operatorname{Hom}(\cdot,\cdot)$ in the second variable, or the left derived functors in the first variable. Properly interpreted, (\ref{extconstruction}) represents only the extreme cases of a general construction \cite{SW}, which attaches $G_{\BB R}$-represen\-tations to all $G_{\BB R}$-orbits. In this way one obtains the maximal globalizations of all standard modules in the Beilinson-Bernstein construction.

Kashiwara \cite{Ka} observed that the results of \cite{SW} could be stated in more functorial language, conjecturally at least, which would then produce not just maximal globalizations of standard modules, but of all Harish-Chandra modules with an infinitesimal character. Proofs of his conjectures appear in \cite{KSch}.

In the remainder of this section we briefly outline the conjectures, respectively the results, of \cite{Ka,KSch}. They involve $\operatorname{\bf D}_{G_{\BB R}} (X)$, the bounded $G_{\BB R}$-equivariant derived category of Bernstein-Lunts \cite{BL}; it is a $G_{\BB R}$-equivariant version of the bounded derived category $\operatorname{\bf D}^b(Sh_X)$ of constructible sheaves \cite{KScha}. Let $S$ be a $G_{\BB R}$-orbit, $j: S \hookrightarrow X$ its inclusion into $X$, and $\mathcal F$ a $G_{\BB R}$-equivariant local system on $S$; then $j_!\mathcal F$, the direct image of $\mathcal F$ with proper supports, is a particular object in $\operatorname{\bf D}_{G_{\BB R}} (X)$. Objects of this type are the basic building blocks, from which the others are put together by successive extensions.

The sheaf of holomorphic sections $\mathcal O(\mathcal L_\lambda)$, $\,\lambda\in\Lambda$, exists as $G$-equivariant sheaf on $X$. More generally, for $\lambda\in\g h^*$ one can make sense of it locally, as an infinitesimally $\g g$-equivariant ``germ of a sheaf" of $\mathcal O$-modules. Collectively these ``germs" constitute a $G$-equivariant {\em twisted sheaf} on $X$, which we denote by $\mathcal O_\lambda$. Technically, $\mathcal O_\lambda$ is not a sheaf on $X$, but rather on the principal $H$-bundle
\begin{equation}
\widehat X\ = \ G/[B_0,B_0]\ \longrightarrow\ G/B_0 \simeq X\,,\ \ \ \text{with structure group}\ \ H = B_0/[B_0,B_0]\,;
\end{equation} 
here $H\simeq (\BB C^*)^r\!$, the {\em universal Cartan group} of $G$, is the connected complex Lie group with Lie algebra $\g h$, and $\widehat X$ is called the {\em enhanced flag variety}. By definition, local sections of $\mathcal O_\lambda$ are locally defined holomorphic functions on $\widehat X$ which transform under right translation by $H$ according to the multiple-valued function $e^\lambda$ on $H$.

Quite analogously one may consider $G_{\BB R}$-equivariant {\em twisted local systems} on $G_{\BB R}$-orbits, with twist $\lambda\in\g h^*$. These, too, are technically sheaves on $\widehat X$, which transform on the right according to the multiple-valued function $e^\lambda$. Whether non-zero twisted local systems, with a particular twist $\lambda$, exist on a particular $G_{\BB R}$-orbit $S$ depends on the interplay between $\lambda$ and the isotropy subgroup $(G_{\BB R})_x\subset G_{\BB R}$ at a point $x\in S$.

Just as the the bounded $G_{\BB R}$-equivariant derived category $\operatorname{\bf D}_{G_{\BB R}} (X)$ is put together from $G_{\BB R}$-equivariant local systems by a process of extensions, the $G_{\BB R}$-equivariant twisted local systems, with twist $\lambda$, are the basic building blocks of the {\em twisted $G_{\BB R}$-equivariant derived category} $\operatorname{\bf D}_{G_{\BB R}} (X)_\lambda$. Let us consider an arbitrary object $\mathcal F\in \operatorname{\bf D}_{G_{\BB R}} (X)_\lambda$. Neither $\mathcal F$ nor $\mathcal O_\lambda$ exist as sheaves on $X$, but both have the {\em same} monodromic behavior under the right action of $H$. Thus homomorphisms between these two twisted sheaves may be regarded as objects on $X$. More precisely, the resolutions from which one would compute $\operatorname{Ext}^*(\mathcal F,\mathcal O_\lambda)$ -- if both arguments were actual sheaves on $X$ -- do exist as complexes of sheaves on $X$. The global $\operatorname{Ext}$ groups
\begin{equation}\label{extconstruction'}
\operatorname{Ext}^p(\mathcal F,\mathcal O_\lambda)\,,\ \ \ \mathcal F\in \operatorname{\bf D}_{G_{\BB R}} (X)_\lambda\,,\ \ \ \lambda\in \g h^*\,,
\end{equation}
are therefore well defined, as complex vector spaces with a linear action of $G_{\BB R}$.
\medskip

\begin{thm}[\cite{KSch}]\label{extconstruction''}
The $\operatorname{Ext}$ groups (\ref{extconstruction'}) carry natural Fr\'echet topologies which make $G_{\BB R}$ act continuously. The resulting representations are admissible, of finite length, and have infinitesimal character $\chi_{\lambda+\rho}$\,. They are the maximal globalizations of their underlying Harish-Chandra modules. The maximal globalization of any Harish-Chandra module with infinitesimal character $\chi_{\lambda+\rho}$ can be realized in this manner.
\end{thm}

For $\mu\in\Lambda$, $e^\mu : H \to \BB C^*$ is well-defined, not multiple-valued. Thus, going back to the definition of the twisted $G_{\BB R}$-equivariant derived category, one finds that $\operatorname{\bf D}_{G_{\BB R}} (X)_\lambda$ depends on $\lambda\in\g h^*$ only modulo the lattice $\Lambda$: 
\begin{equation}\label{periodicity}
\operatorname{\bf D}_{G_{\BB R}} (X)_\lambda\ \ \simeq\ \ \operatorname{\bf D}_{G_{\BB R}} (X)_{\lambda+\mu}\ \ \ \ \text{if}\ \ \mu\in\Lambda\,.
\end{equation}
In particular, for $\lambda\in\Lambda$, $\,\operatorname{\bf D}_{G_{\BB R}} (X)_\lambda \simeq\operatorname{\bf D}_{G_{\BB R}} (X)$. The roles of $\mathcal F$ and $\mathcal O_\lambda$ in (\ref{extconstruction'}) can therefore be played by the constant sheaf $\BB C_X$ and the sheaf of holomorphic sections $\mathcal O(\mathcal L_\lambda)$, with $\lambda\in\Lambda$. In that case the theorem reduces to the Borel-Weil-Bott theorem. Similarly (\ref{extconstruction}) reduces to two special cases of Theorem \ref{extconstruction''}.

\begin{subsection}
{Matsuki Correspondence}
\end{subsection}

If $Q\subset X$ is a $K$-orbit and $S\subset X$ a $G_{\BB R}$-orbit, $\,K_{\BB R}= K\cap G_{\BB R}$ operates on the intersection $Q\cap S$. One calls $Q$ and $S$ {\em dual in the sense of Matsuki} if $Q \cap S$ consists of exactly one $K_{\BB R}$-orbit. The relation ``contained in the closure of" partially orders the set of all $K$-orbits,
\begin{equation}
Q \ \succcurlyeq \ Q'\ \ \ \Longleftrightarrow\ \ \ Q^{\text{cl}}\, \supset \, Q'\,,
\end{equation}
and in the same way orders the set of all $G_{\BB R}$-orbits. 
\bigskip

\begin{thm} [Matsuki \cite{Ma}]
The notion of duality between orbits induces a bijection
$$
\{\text{$G_{\BB R}$-orbits on $X$}\} \ \ \longleftrightarrow\ \ 
\{\text{$K$-orbits on $X$}\}\,,
$$
which reverses the closure relationships.
\end{thm}

As in (\ref{su1,1;1}), let us consider the example of $G_{\BB R}=SU(1,1)$, $G=SL(2,\BB C)$, $K\simeq \BB C^*$. The flag variety of $\g g$ is $X=\BB{CP}^1\simeq \BB C \cup \{\infty\}$, on which $a\in\BB C^*\simeq K$ acts as multiplication by $a^2$. Thus $\{0\}$, $\{\infty\}$, and $\BB C^*$ are the $K$-orbits. The group $SU(1,1)$ acts transitively on the unit disc $D$, the complement $X-D^{\text{cl}}$ of the closure of $D$, and their common boundary $\partial D\simeq S^1$.
\begin{figure}
\centerline{\psfig{figure=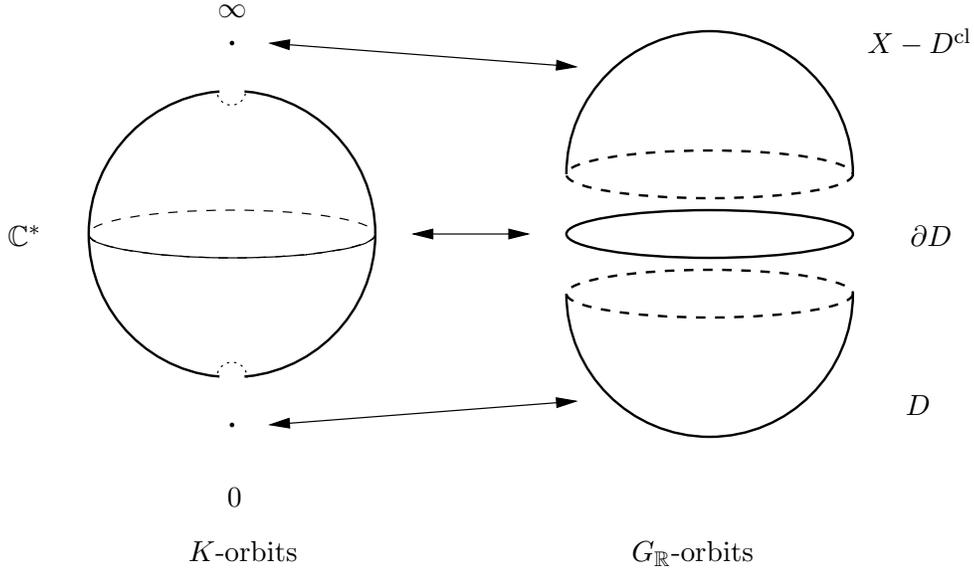}}
\vskip-2.4in \hskip1.8in $\infty$

\hskip5.2in $X-D^{\text{cl}}$

\vskip.82in \hskip.7in ${\BB C}^*$ \hskip4.5in $\partial D$

\vskip.7in \hskip5.4in $D$

\vskip.29in \hskip1.85in $0$

\vskip.1in
\centerline{$K$-orbits \hskip1.7in $G_{\BB R}$-orbits}
\caption{Matsuki correspondence for $SU(1,1)$.}  \label{fig1}
\end{figure}
In this particular situation, duality means that one of the two orbits contains the other, but that is not the case for a general group $G_{\BB R}$.

Matsuki's proof shows that $K$-equivariant twisted local systems on a $K$-orbit $Q$ correspond bijectively to the $G_{\BB R}$-equivariant twisted local systems on the dual $G_{\BB R}$-orbit $S$, in both cases with the same twist $\lambda \in \g h^*$. As was mentioned, the $G_{\BB R}$-equivariant twisted local systems with twist $\lambda$ may be regarded as the basic building blocks of the twisted $G_{\BB R}$-equivariant derived category $\operatorname{\bf D}_{G_{\BB R}} (X)_\lambda$. Quite analogously, the $K$-equivariant derived category $\operatorname{\bf D}_{K} (X)_\lambda$ is built up from $K$-equivariant twisted local systems with twist $\lambda$. But it is not at all obvious that the bijection between equivariant twisted local systems carries over to the extensions between them in the two categories. That was proved by Mirkovi\'c, Uzawa, and Vilonen:
\medskip

\begin{thm} [Matsuki correspondence of sheaves \cite{MUV}]
Matsuki duality induces an equivalence of categories 
$$
\Psi\ : \ \operatorname{\bf D}_{G_{\BB R}} (X)_\lambda \ \ \overset{\sim}{\longrightarrow}\ \ 
\operatorname{\bf D}_K (X)_\lambda\,.
$$
\end{thm}

The covariant deRham functor \cite{Ka1,Me} establishes an equivalence of categories between the categories of, respectively, regular holonomic $\mathcal D$-modules and perverse sheaves constructible with respect to algebraic stratifications. This is the so-called {\em Riemann-Hilbert correspondence}\footnote{The reader may consult Borel's book \cite{Bo} for a discussion of algebraic $\mathcal D$-modules in general and the Riemann-Hilbert correspondence in particular.}. In the context of an algebraic group action with finitely many orbits, all coherent equivariant $\mathcal D$-modules are automatically regular holonomic. The deRham functor therefore restricts to a well defined functor
\begin{equation}\label{deRham}
\operatorname{dR}\ :\ \left\{ \begin{matrix}
\text{category of coherent} \\
\text{$({\cal D},K)$-modules} \end{matrix} \right\}\ \ \longrightarrow\ \ 
\operatorname{\bf D}_K (X)\,.
\end{equation}
For our purposes, it will not matter that (\ref{deRham}) defines an equivalence of categories to the subcategory of perverse objects in $\operatorname{\bf D}_K (X)$. However, we do need the twisted version of the deRham operator,
\begin{equation}\label{twisteddeRham}
\operatorname{dR}\ :\ \left\{ \begin{matrix}
\text{category of coherent} \\
\text{$({\cal D}_{\lambda},K)$-modules} \end{matrix} \right\}\ \ \longrightarrow\ \ 
\operatorname{\bf D}_K (X)_{-\lambda}\,;
\end{equation}
it is contravariant with respect to the twisting, hence the appearance of $\lambda$ on the left and $-\lambda$ on the right. Our final statement relates the Beilinson-Bernstein construction to the construction of $G_{\BB R}$-representations in section 4.4. Like Theorem \ref{extconstruction''}, it was conjectured by Kashiwara \cite{Ka} and proved in \cite{KSch}.
\medskip

\begin{thm}[\cite{KSch}]\label{extconstruction'''}
If $\mathcal S$ is a coherent $({\cal D}_{\lambda},K)$-module, the minimal globalization of the Harish-Chandra module $H^p(X,\mathcal S)$ is isomorphic, as $G_{\BB R}$-representation, to the strong dual of $\,\, \operatorname{Ext}^{n-p}(\Psi^{-1}\circ\operatorname{dR}\,\mathcal S\,,\,\mathcal O_{-\lambda-2\rho})$, with $n=\dim_{\BB C}X$, for any $p\in\BB Z$.
\end{thm}

The sheaf $\,\Psi^{-1}\circ\operatorname{dR}\,\mathcal S\,$ has twist $-\lambda$, but $\,\operatorname{\bf D}_K (X)_{-\lambda} \simeq \operatorname{\bf D}_K (X)_{-\lambda-2\rho}\,$ by the $K$-analogue of (\ref{periodicity}), so it does make sense to consider extensions between $\,\Psi^{-1}\circ\operatorname{dR}\,\mathcal S\,$ and $\,\mathcal O_{-\lambda-2\rho}\,$. 

Some properties of representations are easier to understand in terms of the Beilinson-Bernstein construction, and others easier in terms of the $G_{\BB R}$-construction. Theorems \ref{extconstruction''} and \ref{extconstruction'''} make it possible to ``play off" the two sides against each other. The two theorems also play a crucial role in the proof of the Barbasch-Vogan conjectures \cite{SchV}.

\separate

\bibliographystyle{plain}

\end{document}